\documentclass[11pt]{article}
\usepackage{latexsym}
\usepackage{amsfonts}
\newtheorem{thm}{Theorem}[section]
\newtheorem{cor}[thm]{Corollary}
\newtheorem{lem}[thm]{Lemma}
\newtheorem{pro}[thm]{Proposition}
\newtheorem{defn}[thm]{Definition}

\title{Factorizations of the Thompson-Higman groups, and 
       circuit complexity }

\author{ Jean-Camille Birget
         \thanks{Supported by NSF grant CCR-0310793.}
       }
\date{\today}
\begin{document}
\maketitle

\begin{abstract}
We consider the subgroup ${\it lp}G_{k,1}$ of length preserving elements 
of the Thompson-Higman group $G_{k,1}$ and we show that all elements of 
$G_{k,1}$ have a unique ${\it lp}G_{k,1} \cdot F_{k,1}$ factorization.
This applies to the Thompson-Higman group $T_{k,1}$ as well. We show that
${\it lp}G_{k,1}$ is a ``diagonal'' direct limit of finite symmetric groups,
and that ${\it lp}T_{k,1}$ is a $k^{\infty}$ Pr\"ufer group.
We find an infinite generating set of ${\it lp}G_{k,1}$ which is related to
reversible boolean circuits. 

We further investigate connections between the Thompson-Higman groups, 
circuits, and complexity.
We show that elements of $F_{k,1}$ cannot be one-way functions.
We show that describing an element of $G_{k,1}$ by a generalized bijective 
circuit is equivalent to describing the element by a word over a certain 
infinite generating set of $G_{k,1}$; word length over these generators is 
equivalent to generalized bijective circuit size.

We give some coNP-completeness results for $G_{k,1}$ (e.g., the word problem
when elements are given by circuits), and $\#{\mathcal P}$-completeness 
results (e.g., finding the ${\it lp}G_{k,1} \cdot F_{k,1}$ factorization 
of an element of $G_{k,1}$ given by a circuit).  
\end{abstract}


\section{Introduction}

The Thompson groups, introduced by Richard J.\ Thompson \cite{Th0,McKTh,Th},
and their generalization by Graham Higman \cite{Hig74}, are well known 
for their amazing properties and their importance in combinatorial group 
theory and topology. In this paper we focus on the computational role of 
these groups, continuing the work started in \cite{BiThomps, BiCoNP}, and
we study some subgroups of the Thompson-Higman groups that are motivated by
circuit complexity.
We emphasize that we view the Thompson groups as {\em a model of computation} 
and not just a source of algorithmic problems. 
Indeed, since the Thompson group $V$ has a faithful partial action 
on the set $\{0,1\}^*$ of all bitstrings, it is natural to consider 
combinational circuits for computing elements of $V$; i.e., we can view 
every element of $V$ as the input-output function of an acyclic digital 
circuit (see \cite{BiCoNP}). The elements of $V$ are bijections, hence we
will see connections between $V$ and {\it reversible computing}. More
precisely, words over certain generating sets of $V$ will be seen to be
equivalent to circuits made of (generalized) bijective gates. 

Combinational circuits have fixed-length inputs and fixed-length outputs, 
which is not the case for elements of $V$; but the notion of a circuit can 
be adapted in order to be applied to the computation of elements of $V$. 
Moreover, for bijective circuits the fixed length of inputs and outputs 
implies that the circuit is length preserving. 
This leads to the question what the {\it length preserving} elements of 
$V$ are, and how arbitrary elements of $V$ are related to length preserving 
elements of $V$. The length preserving elements of $V$ turn out to form an
interesting subgroup, called ${\it lp}V$, and we will show that every 
element of $V$ can be factored in a unique way as a product of an element 
of ${\it lp}V$ and an element of the Thompson group $F$. This factorization 
carries over to the Thompson group $T$, where we have a unique 
${\it lp}T \cdot F$ factorization. All of this generalizes to the 
Thompson-Higman groups $G_{k,1}$ and $T_{k,1}$.  

The group ${\it lp}G_{k,1}$ is locally finite; it is a ``diagonal'' direct 
limit of finite symmetric groups, and it is simple when $k$ is even; 
${\it lp}T_{k,1}$ is a $k^{\infty}$ Pr\"ufer group. 
The connection with bijective (a.k.a.\ ``reversible'') circuits leads to an
interesting infinite generating set of ${\it lp}V$. 
We show that a description of an element of $G_{k,1}$ by a bijective circuit 
is equivalent to a description by a word over a certain infinite generating 
set of $G_{k,1}$; bijective circuit size is closely related to the word size
over a certain infinite generating set of $G_{k,1}$. This shows that $G_{k,1}$
(and especially $V$) can serve as a model for bijective computing, with 
equivalent complexity.  

We also investigate the computational complexity of some problems in
$G_{k,1}$. We show that when an element $\varphi \in G_{k,1}$ is given by a 
bijective circuit (or by a general non-bijective circuit), the question 
whether $\varphi$ is the identity, and the question whether $\varphi$ is 
maximally extended, are coNP-complete problems; this is an application of 
\cite{BiCoNP} where $G_{k,1}$ and its connection with circuits was used to 
construct a finitely presented group with coNP-complete word problem.  
We show that elements of $F_{k,1}$ cannot be one-way functions, i.e., from 
a circuit for $f \in F_{k,1}$ one can easily find a circuit for $f^{-1}$.
And we show that when $\varphi \in G_{k,1}$ is given by a bijective circuit,
the problem of finding the ${\it lp}G_{k,1} \cdot F_{k,1}$ factorization of
$\varphi$ is $\#{\mathcal P}$-complete in general (and also under some
restrictions). 

\bigskip

\noindent {\bf Definition of the Thompson-Higman group} 

The rest of this Introduction consists of a brief, but complete, definition 
of the Thompson-Higman groups $G_{k,1}$, $T_{k,1}$, and $F_{k,1}$. 
We follow the exposition of \cite{BiThomps,BiCoNP}, based on partial actions 
on finite words, which simplifies the connections with circuits. Compare with 
the definition in \cite{Th} (based on infinite sequences), \cite{Hig74} 
(based on automorphisms of certain algebras), \cite{Scott} (based on words 
and similar to this paper, but with different terminology), \cite{CFP} 
(based on finite trees), \cite{BrinSqu} (based on piecewise linear maps 
between real numbers, see also \cite{CFP}), \cite{Dehornoy} (related to 
associativity or commutativity in term rewriting, which was Thompson's 
original view). 

To define the Thompson-Higman group $G_{k,1}$ we fix an alphabet $A$ of 
cardinality $|A| = k$. Let $A^*$ denote the set of all finite {\it words} 
over $A$ (i.e., all finite sequences of elements of $A$); this includes the 
empty word $\varepsilon$. The {\it length} of $w \in A^*$ is 
denoted by $|w|$; let $A^n$ denote the set of words of length $n$. For two
words $u,v \in A^*$ we denote their {\it concatenation} by $uv$ or by 
$u \cdot v$; for sets $B, C \subseteq A^*$ the concatenation is 
$BC = \{uv : u \in B, v \in C\}$.
A {\it right ideal} of $A^*$ is a subset $R \subseteq A^*$ such that 
$RA^* \subseteq R$. A generating set of a right ideal $R$ 
is a set $C$ such that $R$ is the intersection of all right ideals that
contain $C$. A right ideal $R$ is called {\it essential} iff $R$ has a 
non-empty intersections with every right ideal of $A^*$. 
For words $u,v \in A^*$, we say that $u$ is a {\it prefix} of $v$ iff there 
exists $z \in A^*$ such that $uz = v$. A {\it prefix code} is a subset 
$C \subseteq A^*$ such that no element of $C$ is a prefix of another element 
of $C$. A prefix code is {\it maximal} iff it is not a strict subset of 
another prefix code.
One can prove that a right ideal $R$ has a unique minimal (under inclusion)
generating set, and that this minimal generating set is a prefix code; this 
prefix code is maximal iff $R$ is an essential right ideal.

Partial functions on $A^*$ will play a big role. For $f: A^* \to A^*$, let 
Dom$(f)$ denote the domain and let Im$(f)$ denote the image (range) of $f$.
A {\it restriction} of $f$ is any function $f_1: A^* \to A^*$ such that 
Dom$(f_1) \subseteq$ Dom$(f)$, and such that $f_1(x) = f(x)$ for all 
$x \in {\rm Dom}(f_1)$. An {\it extension} of $f$ is any function on $A^*$ 
of which $f$ is a restriction. 

An {\it isomorphism} between right ideals $R_1, R_2$ of $A^*$ is a bijection
$\varphi: R_1 \to R_2$ such that for all $r_1 \in R_1$ and all $z \in A^*$:
 \ $\varphi(r_1z) = \varphi(r_1) \cdot z$; the isomorphism $\varphi$ can be
described by a bijection between the prefix codes that minimally generate
$R_1$, respectively $R_2$. 

One can prove that an isomorphism $\varphi$ between essential right ideals 
has a unique maximal extension (as an isomorphism between essential right 
ideals), denoted max $\varphi$. So, max $\varphi$ has no extension (other 
than itself) to an isomorphism between essential right ideal. 

Finally, the Thompson-Higman group $G_{k,1}$ is defined to consist of all 
maximally extended isomorphisms between finitely generated essential right 
ideals of $A^*$.
The multiplication consists of composition followed by maximal extension:
$\varphi \cdot \psi = $ max$(\varphi \circ \psi)$. 
Note that we let $G_{k,1}$ act partially and faithfully on $A^*$ on the 
{\it left}. 

Thompson and Higman proved that $G_{k,1}$ is finitely presented. Also, when 
$k$ is even $G_{k,1}$ is simple, and when $k$ is odd $G_{k,1}$ has a simple 
normal subgroup of index 2. 

Every element $\varphi \in G_{k,1}$ can be described by a bijection between 
two finite maximal prefix codes; this bijection can be described concretely 
by a finite function {\it table}. When $\varphi$ is described 
by a maximally extended isomorphism between essential right ideals, 
$\varphi: R_1 \to R_2$, we call the minimum generating set of $R_1$ the 
{\it domain code} of $\varphi$, denoted domC$(\varphi)$, and we call the
minimum generating set of $R_2$ the {\it image code} of $\varphi$, denoted 
imC$(\varphi)$; because of the uniqueness of maximal extension,
domC$(\varphi)$ and imC$(\varphi)$ are uniquely determined by $\varphi$. 
We call the cardinality $|{\rm domC}(\varphi)| = $
$|{\rm imC}(\varphi)|$ the {\it table size} of $\varphi$, denoted 
$\|\varphi\|$. In \cite{BiThomps} it was proved that for all 
$\varphi, \psi \in G_{k,1}$: \ \ 
$\|\varphi \psi\| \leq \|\varphi\| + \|\psi\|$.  The concepts of
domC$(\varphi)$, imC$(\varphi)$, table, and $\|\varphi\|$, can also be used 
when $\varphi$ is not maximally extended.

For any finite generating set $\Gamma$ of $G_{k,1}$ and any 
$\varphi \in G_{k,1}$, we define the {\it word length} of $\varphi$ over 
$\Gamma$ as the length of a shortest word over $\Gamma \cup \Gamma^{-1}$ that 
represents $\varphi$; it is denoted by $|\varphi|_{\Gamma}$. In 
\cite{BiThomps} it was proved that for any 
finite generating set $\Gamma$ of $G_{k,1}$, the word length and the 
table size are closely related; for all $\varphi \in G_{k,1}$: \ \
$c' \ \|\varphi\| \ \leq \ |\varphi|_{\Gamma} \ \leq \ $
$c \ \|\varphi\| \ \log_2 \|\varphi\|$ \ (for some constants $c, c' >0$ 
depending on $\Gamma$ but not on $\varphi$). Asymptotically, for most 
$\varphi \in G_{k,1}$ we also have \ $|\varphi|_{\Gamma} \ \geq \ $
$c'' \ \|\varphi\| \ \log_2 \|\varphi\|$ \ (for some constant depending on 
$\Gamma$, $0 < c'' < c$).  However, for $\varphi \in F_{k,1}$ it was proved 
in \cite{CFP} that \ $c' \ \|\varphi\| \ \leq \ |\varphi|_{\Gamma} \ $
$ \leq \ c \ \|\varphi\|$. 

\medskip

We will use the well-known finitely presented subgroups $F_{k,1}$ and
$T_{k,1}$ of $G_{k,1}$, introduced in \cite{Th0} and \cite{Hig74}. The 
groups $F_{2,1}$ (also called $F$) and $T_{2,1}$ (also called $T$)
have a large literature; a few examples are \cite{McKTh}, \cite{McKTh,Th}, 
\cite{BrinSqu}, \cite{CFP}, \cite{GhysSergiescu}, \cite{BrownGeo}, 
\cite{ClearyTaback}, \cite{Brin97}, \cite{Brin99}, \cite{GubaSapir}, 
\cite{BCST}. 
Below we will introduce the subgroups ${\it lp}G_{k,1}$ and ${\it lp}T_{k,1}$
of length preserving elements of $G_{k,1}$, respectively $T_{k,1}$. 

We will need the exact definition of $F_{k,1}$ and $T_{k,1}$ in the setting 
of partial actions on words (in $A^*$), and to do so we need some preliminary
definitions. Assuming that a linear order has been chosen for the alphabet 
$A$, we can consider the dictionary order on $A^*$, denoted $\leq_{\rm d}$,
and defined as follows. For any 
$x_1, x_2 \in A^*$ we say that $x_1 \leq_{\rm d} x_2$ (i.e., $x_1$ precedes
$x_2$ in the dictionary order) iff either (1) \ $x_1$ is a prefix of $x_2$, 
or, (2) \ letting $p$ denote the longest common prefix of $x_1$ and $x_2$, 
we have: \ $x_1 = pa_1v_1$, \ $x_2 = pa_2v_2$, \ with \ $a_1 < a_2$ 
(for some letters $a_1, a_2 \in A$, and words $v_1, v_2 \in A^*$, where $<$ is 
the strict order in $A$).

A partial map $f: A^* \to A^*$ is said to {\it preserve the 
dictionary order} iff for all $x_1, x_2 \in {\rm Dom}(f)$ we have: \   
$x_1 \leq_{\rm d} x_2$ \ iff \ $f(x_1) \leq_{\rm d} f(x_2)$. 

We also want to define ``cyclical preservation'' of the dictionary order. 
Here we will simply write $<$ for $<_{\rm d}$ (strict dictionary order). 
A cyclical order of a finite maximal prefix code $P \subset A^*$ is a listing
$(x_0, x_1, \ldots, x_{|P|-1})$ of all the elements of $P$ such that for some
integer $s$: \ $(x_s, \ldots, x_{|P|-1}, x_0, \ldots, x_{s-1})$ is the 
listing of $P$ in dictionary order. In other words, a cyclical order of $P$ 
is a cyclic permutation of the dictionary order on $P$. 

We say that a partial map $f: A^* \to A^*$ {\it cyclically preserves the 
dictionary order} iff for all finite sequences $(x_0, x_1,  \ldots, x_{n-1})$
we have: \ $(x_0, x_1,  \ldots, x_{n-1})$ is a cyclical order of some finite 
maximal prefix code iff $(f(x_0), f(x_1), \ldots, f(x_{n-1}))$ is a cyclical 
order of some finite maximal prefix code.  

The groups $F_{k,1}$ and $T_{k,1}$ can be defined as follows, from the
point of view of partial actions on finite words (see \cite{BiThomps}).

\begin{defn} \label{F_and_T}
Assume that a linear order has been chosen for the alphabet $A$, where 
$|A| = k$. Then $F_{k,1}$ consists of the elements of $G_{k,1}$ that 
preserve the dictionary order of $A^*$, and 
$T_{k,1}$ consists of the elements of $G_{k,1}$ that cyclically preserve
the dictionary order of $A^*$. 
\end{defn}

\noindent {\bf Another view of $F_{k,1}$:} The elements of $F_{k,1}$ can be 
given the following interpretation. First we define the concept of a rank
function on a (partial) order structure $(S,\leq)$. The {\em rank} of an 
element $t \in S$ is \ 

\smallskip

 \ \ \ \ \  ${\rm rank}_S(t) = |\{ x \in S : x < t \}|$, 

\smallskip

\noindent i.e., the number of elements that strictly precede $t$.
Every element of $\varphi \in G_{k,1}$ can be represented (after appropriate 
restriction) by a bijective partial function $\varphi: A^* \to A^*$ such that 
${\rm imC}(\varphi) = A^n$ for some $n>0$, and ${\rm domC}(\varphi)$ is some 
finite maximal prefix code of cardinality $k^n$ (where $|A| = k$). 
If we view the elements of $A^n$ as the integers $\{0,1, \ldots, k^n -1 \}$ 
in base-$k$ representation we have: 

\smallskip

 \ \ \ $F_{k,1}$ {\it consists of all elements of $G_{k,1}$ that can be
 represented by rank functions}

 \ \ \ ${\rm rank}_P(.): P \to \{0,1, \ldots, k^n -1 \}$, 

\smallskip

\noindent where $n$ ranges over the positive integers, and
$P \subset A^*$ ranges over all maximal prefix codes of cardinality $k^n$.
This point of view will help us later in proving that elements of $F_{k,1}$
can have high computational complexity, even when their domain code 
${\rm domC}(\varphi)$ has an easy membership problem (see Theorem 
\ref{numberP_fact2}).


\bigskip

\noindent {\bf Overview:} 
This paper consists of the following parts: 

Part 1 consists of sections 2, 3, and 4. We introduce the subgroup
${\it lp}G_{k,1}$ of length preserving elements of the Thompson-Higman group
$G_{k,1}$, and we give the ${\it lp}G_{k,1} \cdot F_{k,1}$ factorization of
$G_{k,1}$; we generalize this unique factorization to other subgroups of
$G_{k,1}$.

Section 5 makes the transition from part 1 to part 2, by giving a connection
between circuits and some properties of ${\it lp}G_{k,1}$.

Part 2 consists of sections 6 and 7.  We study $V$ as a model for reversible 
circuits. We also investigate the complexity of some problems: We show that 
elements of $F$ cannot be one-way functions, and we show that finding the 
${\it lp}V \cdot F$ factorization of an element of $V$ given by a circuit is
$\#{\mathcal P}$-complete.


\section{The subgroups ${\it lp}G_{k,1}$ and ${\it lp}T_{k,1}$ }

The Thompson-Higman group $G_{k,1}$ contains all finite symmetric groups, and 
this inspires the definition of the subgroup ${\it lp}G_{k,1}$ of all 
length-preserving elements of $G_{k,1}$. We will denote ${\it lp}G_{2,1}$ 
also by ${\it lp}V$.
Another motivation of ${\it lp}V$, which we will develop more later, is the 
computation of elements of $V$ and ${\it lp}V$ by digital circuits. Indeed, 
circuits traditionally have a fixed length for inputs and a fixed length for 
outputs (corresponding to fixed numbers of wires); for bijective functions 
this means length preservation. 

\begin{defn}
The subgroup of length-preserving elements of the Thompson-Higman group
$G_{k,1}$ is \ lp$G_{k,1} = \{ \varphi \in G_{k,1} : $
$\forall x \in$ {\rm Dom}$(\varphi), \, |x| = |\varphi(x)| \}$.
 \ Similarly we define \ lp$T_{k,1} = T_{k,1} \, \cap \, {\it lp}G_{k,1}$.
\end{defn}
Restriction or extension of a length-preserving partial function
$A^* \to A^*$, representing an element of $G_{k,1}$, is again length
preserving, so {\it lp}$G_{k,1}$ is well-defined as a subset of the group 
$G_{k,1}$. The inverse of a length-preserving partial function  is also 
length-preserving.
After a restriction, if necessary, any finite set of elements of 
{\it lp}$G_{k,1}$ can be represented by permutations of the same set $A^m$, 
for any large enough $m$.
Hence {\it lp}$G_{k,1}$ is closed under composition. It follows that  
{\it lp}$G_{k,1}$ is a subgroup of $G_{k,1}$.  
The group {\it lp}$G_{k,1}$ is {\it locally finite} (i.e., every finitely 
generated subgroup is finite), and {\it lp}$G_{k,1}$ contains all the finite 
symmetric groups ${\mathfrak S}_{A^n}$, for all $n \geq 1$.


\bigskip

Assume we restrict an element $\varphi \in {\it lp}G_{k,1}$ so that its 
domain and image codes are both $A^m$ for some $m$. 
Then the additional overall restriction operation (which replaces each 
$\varphi(x) = y$ by the $k$-tuple $\varphi(x \, a) = y \, a$, where 
$a$ ranges over $A$) leads to the following embeddings:

\smallskip

 \ \ \ \ \ 
$\otimes {\rm id}_A: \pi \in {\mathfrak S}_{A^n} \ \hookrightarrow \ $
$\pi \otimes {\rm id}_A \in {\mathfrak S}_{A^{n+1}}$, \  

\smallskip

\noindent where for all 
$x \in A^n$ and $a \in A$ we define \   
$(\pi \otimes {\rm id}_A)(x a) = \pi(x) \cdot a$ (where $\cdot$ denotes 
concatenation). This type of embedding of symmetric groups is called 
{\it diagonal} \cite{Zalesskii}, \cite{Hartley}.  
Moreover, when $|A| = k$ is {\em even} then the above embedding factors 
through the alternating group \

\smallskip

 \ \ \ \ \  ${\mathfrak S}_{A^n} \ \hookrightarrow \ $
${\mathfrak A}_{A^{n+1}} \ \subset {\mathfrak S}_{A^{n+1}}$.

\medskip

\noindent Indeed we have the following generalization of an observation of
\cite{Shende} (see also Section 5 below): For any positive integer $n$ and 
any $\pi \in {\mathfrak S}_{A^n}$, the permutation $\pi \otimes {\rm id}_A$ 
is {\it even}. 
Indeed, the transformation $\pi \to \pi \otimes {\rm id}_A$  replaces one 
transposition $(u|v)$ of $\pi$ (with $u,v \in A^n$) by the sequence of 
transpositions $(ua_1|va_1) \ldots (ua_k|va_k)$,  i.e., $k$ 
transpositions with $k$ even; here $A = \{a_1, \ldots, a_k\}$.

\medskip

\noindent 
The above embeddings yield the following. 

\begin{pro} \label{lpV_infinite_alt} 
The group ${\it lp}G_{k,1}$ is isomorphic to the direct limit of the 
sequence of diagonal embeddings \  
$\otimes {\rm id}_A: {\mathfrak S}_{A^n} \ \hookrightarrow$ 
$ \ {\mathfrak S}_{A^{n+1}}$. 

When $|A| = k$ is even, ${\it lp}G_{k,1}$ is isomorphic to the direct limit of
the sequence of embeddings \ 
${\mathfrak S}_{A^n} \ \hookrightarrow \ $
${\mathfrak A}_{A^{n+1}} \ \subset {\mathfrak S}_{A^{n+1}}$.
\end{pro}

These are examples of the direct limits of symmetric groups considered in 
\cite{KegelWehrfritz}, chapter 6, and in \cite{Hartley}, section 1.5. 
The embedding maps are of ``diagonal'' type, in the terminology 
of these references. By these references we also conclude that when $k$ is 
even, ${\it lp}G_{k,1}$ is a simple group, and when $k$ is odd, 
${\it lp}G_{k,1}$ has a simple subgroup of index 2 (via the parity map).
In any case, it also follows from \cite{KegelWehrfritz} and \cite{Hartley} 
that ${\it lp}G_{k,1}$ is different from the finitary symmetric group and the 
finitary alternating group; indeed, the finitary symmetric group does not 
contain any Pr\"ufer groups, whereas ${\it lp}G_{k,1}$ contains 
${\it lp}T_{k,1}$ which is a group of Pr\"ufer type (as we shall see next).
However, ${\it lp}G_{k,1}$ also contains many copies of the finitary 
symmetric group (as was mentioned in \cite{Th}). 
 
\bigskip

The observations above apply also to the Thompson group $T_{k,1}$.
Let us denote by ${\mathbb Z}_{A^n}$ the cyclic subgroup of 
${\mathfrak S}_{A^n}$ generated by the permutation  \   
$w_i \mapsto w_{(i+1) \, {\rm mod} \, k^n}$, where 
$(w_i : i = 0,1, \ldots, k^n -1)$ is the listing of $A^n$ in dictionary order.  
${\mathbb Z}_{A^n}$ consists of the elements of ${\mathfrak S}_{A^n}$
that cyclically preserve the dictionary order.  
Just as for the symmetric groups on $A^n$, the restriction operation of 
$G_{k,1}$ gives an embedding of ${\mathbb Z}_{A^n}$ into 
${\mathbb Z}_{A^{n+1}}$, by the transformation $\otimes {\rm id}_A$ which 
sends the generator 
 \ $(w_i \mapsto w_{(i+1) \, {\rm mod} \, k^n})$ \ of ${\mathbb Z}_{A^n}$ to 
the element $(v_j \mapsto v_{(j+k) \, {\rm mod} \, k^{n+1}})$ of 
${\mathbb Z}_{A^{n+1}}$.  Here, $(v_j : j = 0,1, \ldots, k^{n+1} -1)$ is the 
listing of $A^{n+1}$ in dictionary order.  Thus we have:

\begin{pro} \label{lpT_Prufer}
The group ${\it lp}T_{k,1}$ is isomorphic to the $k^{\infty}$ Pr\"ufer group,
given by the direct limit of the sequence of embeddings \   
${\mathbb Z}_{A^n} \hookrightarrow {\mathbb Z}_{A^{n+1}}$, 
where the embeddings are determined by the restriction operation of $G_{k,1}$.
\end{pro}

The $k^{\infty}$ Pr\"ufer group is isomorphic to the multiplicative group of
the complex $k^n$th roots of unity (for all $n>0$), or the additive group of
$k$-ary rationals modulo 1, i.e.,
$\{ \frac{m}{k^n} \ {\rm mod} \ 1 \ : \ n, m \in {\mathbb N} \}$.

\section{Length-preserving order-preserving factorization of $G_{k,1}$ 
          and $T_{k,1}$ }

Let {\bf 1} denote the identity of $G_{k,1}$.
\begin{lem} \label{lp_cap_F}
If an element of $F_{k,1}$ has a representation $f: A^* \to A^*$ such that
${\rm domC}(f) = {\rm imC}(f)$ then $f$ represents the identity.
Hence, \
  $F_{k,1} \ \cap \ {\it lp}G_{k,1} \ = \ \{ {\bf 1} \}.$
\end{lem}
{\bf Proof.} \ Since ${\rm domC}(f) = {\rm imC}(f)$, $f$ is a permutation of
${\rm domC}(f)$. For any finite set of words, the only permutation that
preserves the dictionary order is the identity.

We saw already that every element $\varphi \in {\it lp}G_{k,1}$ can be
represented by a permutation of $A^m$ for some $m>0$, so
${\rm domC}(\varphi) = {\rm imC}(\varphi)$ for every 
$\varphi \in {\it lp}G_{k,1}$.     \ \ \ $\Box$

\begin{thm} \label{factorization_lpV_F}. \\     
$\bullet$ We have \ 
 $G_{k,1} \ = \ {\it lp}G_{k,1} \cdot F_{k,1}$
 \ where every element $\varphi$ of $G_{k,1}$ has a {\em unique}
factorization \ $\varphi = \pi \cdot f$ with $\pi \in {\it lp}G_{k,1}$ and
$f \in F_{k,1}$.

\noindent $\bullet$ Symmetrically there is a unique factorization \
$G_{k,1} \ = \  F_{k,1} \cdot {\it lp}G_{k,1}$.

\noindent $\bullet$ For $T_{k,1}$ there are unique factorizations \ 
$T_{k,1} \ = \ {\it lp}T_{k,1} \cdot F_{k,1} $ $ \ = \ $ 
$F_{k,1} \cdot {\it lp}T_{k,1}$.
\end{thm}
{\bf Proof.} \ Uniqueness of the factorization follows immediately from
Lemma \ref{lp_cap_F}: If $\pi_1 f_1 = \pi_2 f_2$ then
$\pi_2^{-1} \pi_1 = f_2 f_1^{-1} \in \, F_{k,1} \cap {\it lp}G_{k,1} = $
$\{ {\bf 1} \}$, \ hence $\pi_2^{-1} \pi_1 = {\bf 1} = f_2 f_1^{-1}$.
Existence follows from the following factorization algorithm, whose input
is any $\varphi \in G_{k,1}$.

\smallskip

\noindent {\it Factorization algorithm}:

\smallskip

\noindent (1) \ \ \ Restrict $\varphi$ so that its image code becomes $A^n$ 
for some $n>0$. Let $P$ be the corresponding domain code (of cardinality 
$k^n$). So now $\varphi$ is represented by a bijection $P \to A^n$.

\smallskip

\noindent (2) \ \ \ Let $f: P \to A^n$ be the unique element of $F_{k,1}$ 
determined by the finite maximal prefix codes $P$ and $A^n$. 

\smallskip

\noindent (3) \ \ \ Let $\pi(.) = \varphi \ f^{-1}(.)$; then 
$\varphi = \pi f$. 
 \ [End of algorithm.]  
 
\smallskip 

We claim that $\pi \in {\it lp}G_{k,1}$.  Indeed, the domain code and the 
image code of $\pi$ are both $A^n$; hence $\pi$ preserves length.

\smallskip

In the case of $T_{k,1}$ we observe that if $\varphi \in T_{k,1}$ then
the unique factorization $\varphi = \pi \, f$ yields 
$\pi = \varphi \, f^{-1} \in T_{k,1}$ (since $\varphi \in T_{k,1}$ and
$F_{k,1} \subset T_{k,1}$).
 \ \ \ $\Box$

\smallskip 

Observe that $f \in F_{k,1}$, produced by the factorization algorithm, is 
the ranking function of $P$, when we view $A^n$ as the integers
$\{0,1, \ldots, k^n -1\}$ in base-$k$ notation.

\bigskip

We will examine how the table sizes of $\pi \in {\it lp}G_{k,1}$ and 
$f \in F_{k,1}$ are related to the table size of $\varphi$ when 
$\varphi = \pi \, f$. 
It turns out that $\pi$ and $f$ can have exponentially larger 
size than $\varphi$. 
In a later section we'll consider other complexity measures for
$\pi$ and $f$.

\begin{pro} \label{table_size_pi_f}
For all $n > 2$ there are elements $\varphi_n \in T_{2,1}$ whose 
factorization $\varphi_n = \pi_n f_n$ leads to an exponential 
increase in table size. More precisely, $\varphi_n$ can be found so that 
$\|\varphi_n\| = n$, and $\|\pi_n\| = \|f_n\| = 2^{n-1}$. 
\end{pro}
{\bf Proof.} \ Let us pick $\varphi_n \in T_{2,1}$ given by the following
table, over the alphabet $A = \{a,b\}$:

\[ \varphi_n \ = \ \left[
\begin{array}{ccc cc ccc}
a^{n-1} \ & \ a^{n-2}b \ & \ \ldots \ & \ a^ib \ & \ \ldots \ & \ a^2b \ & \
ab \ & \ b  \\
a^{n-2}b \ & \ a^{n-3}b \ & \ \ldots \ & \ a^{i-1}b \ & \ \ldots \ & \ ab \ &
 \ b \  & \ a^{n-1}
\end{array}        \right]
\]

\noindent So, $\varphi_n$ is a cyclic permutation of the finite maximal 
prefix code \   
$\{a^{n-1}\} \cup \{a^i b : i = n-2, \ldots , 1, 0 \}$. One observes that
$\varphi_n$ is reduced (unextendable) as given by the table, hence
$\|\varphi_n\| = n$.

The longest words in the image code of $\varphi_n$ in the above table
have length $n-1$. When we restrict $\varphi_n$ and let its image code 
become $\{a,b\}^{n-1}$ we obtain the following table of size $2^{n-1}$ for
$\varphi_n$, where $x_j$ ranges over $\{a,b\}^j$ (for $j=1, \ldots, n-2$):

\[ \varphi_n \ = \ \left[
\begin{array}{ccc cc ccc}
a^{n-1} \ & \ a^{n-2}b x_1 \ & \ \ldots \ & \ a^ib x_{n-i-1} \ & \ \ldots \
& \ a^2b x_{n-3} \ & \ ab x_{n-2} \ & \ b  \\
a^{n-2}b \ & \ a^{n-3}b x_1 \ & \ \ldots \ & \ a^{i-1}b x_{n-i-1} \ & \ \ldots
 \ & \ ab x_{n-3} \  & \ b x_{n-2} \  & \ a^{n-1}
\end{array}        \right]
\]

\bigskip

\noindent Then in the factorization $\varphi_n(.) = \pi_n f_n(.)$ we have:

\medskip

\[ f_n \ = \            \left[
\begin{array}{ccc ccc cc}
a^{n-1} \ & \ a^{n-2}b \, a \ & \  a^{n-2}b \, b \ & \ \ldots
 \ & \ a^ib \, a^{n-i-1} \ & \ a^ib \, s(x_{n-i-1}) \ & \ a^{i-1}b a^{n-i}
   \ & \ \ldots  \\
a^{n-1} \ & \ a^{n-2}b \      & \ a^{n-3}b \, a \  & \ \ldots
  \ & \ a^ib \, b^{n-i-2} \ & \ a^{i-1}b x_{n-i-1} \ & \ a^{i-1}b \, b^{n-i-1}
   \ & \ \ldots \
\end{array}    \right.
\]

\medskip

\hspace{3in}
\[            \left.
\begin{array}{ccc ccc}
 \ \ldots \ & \ a^2b \, a^{n-3} \ & \ a^2b \, s(x_{n-3})
      \ & \ ab \, a^{n-2} \ & \ ab \, s(x_{n-2}) \ & \ b  \\
 \ \ldots  \ & \ a^2b \, b^{n-4} \ & \ ab x_{n-3}
      \ & \ b \, b^{n-3}  \ & \ bx_{n-2} \  & \ b^{n-1}
\end{array}        \right]
\]

\medskip

\noindent where each $x_j$ ranges over $\{a,b\}^j - \{b^j\}$, and where
$s(x_j)$ denotes the {\it successor} of $x_j$ in the dictionary order; hence,
$s(x_j)$ ranges over $\{a,b\}^j - \{a^j\}$. In the table, the strings $x_j$
and the strings $s(x_j)$ appear in dictionary order. We also have

\medskip

\medskip

\[ \pi_n \ = \            \left[
\begin{array}{ccc ccc cc}
b^{n-1} \ & \ a^{n-1}  \ & \ a^{n-2}b \      & \ a^{n-3}b \, a \ & \ \ldots
 \ & \ a^ib \, b^{n-i-2} \     & \ a^{i-1}b \, x_{n-i-1} \ & \ \ldots  \\
a^{n-1} \ & \ a^{n-2}b \  & \ a^{n-3}b \, a \ & \ a^{n-3}b \, b \ & \ \ldots
  \ & \ a^{i-1}b \, a^{n-i-1} \ & \ a^{i-1}b \, s(x_{n-i-1}) \ & \ \ldots
\end{array}    \right.
\]

\medskip

\hspace{3in}
\[            \left.
\begin{array}{ccc ccc}
 \ldots \ & \ a^2b \, b^{n-4} \ & \ a^2b \, x_{n-3}
      \ & \ ab \, b^{n-3} \ & \ b \, x_{n-2}  \\
 \ldots \ & \ ab \, a^{n-3} \ & \ ab \, s(x_{n-3})
      \ & \ b \, a^{n-2}  \ & \ b \, s(x_{n-2})
\end{array}        \right]
\]

\medskip

\noindent where the words $x_j$ and $s(x_j)$ range over the same values and 
have the same meaning as for $f_n$.

\smallskip

One sees in the table of $\pi_n$ that for every argument $x$, $\pi_n(x)$ 
differs from $x$ in the right-most letter: whenever $x$ ends in $a$, 
$\pi_n(x)$ ends in $b$, and vice versa. Hence, $\pi_n$ as given by the table, 
is reduced (cannot be extended). Hence, $\|\pi_n\|$ is the size of the above
table, i.e., $2^{n-1}$.
Similarly, in the table for $f_n$, $x$ and $f_n(x)$ differ in the right-most
letter, except when $x = a^{n-1}$ or $x = b$. Hence $f_n$ as given by the
table is reduced, and $\|f_n\| = 2^{n-1}$.
 \ \ \ $\Box$


\section{Other factorizations of $G_{k,1}$ and $T_{k,1}$}

We will give an infinite collection of torsion subgroups $S$ of $G_{k,1}$ 
that can be used for factoring $G_{k,1}$ as $S \cdot F_{k,1}$.

If $P \subset A^*$ is a finite maximal prefix code then for every $n \geq 0$, 
the overall restriction operation in $G_{k,1}$ determines a diagonal 
embedding 

\smallskip 

  \ \ \ $\otimes {\rm id}_A: \pi \in {\mathfrak S}_{PA^n} \ \hookrightarrow$
 $ \ \pi \otimes {\rm id}_A \in {\mathfrak S}_{PA^{n+1}}$

\smallskip

\noindent where $(\pi \otimes {\rm id}_A)(xa) = \pi(x) \cdot a$, for all 
$x \in PA^n$, $a \in A$. This is a generalization of the embedding  
${\mathfrak S}_{A^n}$ $\hookrightarrow {\mathfrak S}_{A^{n+1}}$ that we saw
earlier (which was the special case when $P$ consists of just the empty word). 
We then take the direct limit of this sequence of symmetric groups 
and obtain a subgroup of $G_{k,1}$, denoted by

\smallskip 

  \ \ \ $\bigcup_{n \geq 1} {\mathfrak S}_{P A^n}$.    

\smallskip

\noindent Just as for ${\it lp}G_{k,1}$, when $k = |A|$ is even the group 
$\bigcup_{n \geq 1} {\mathfrak S}_{P A^n}$ is simple, and when $k$ is odd
the group has a simple subgroup of index 2 (via the parity map).

\begin{thm} \label{PAm_factorization}
If $P \subset A^*$ is a finite maximal prefix code then 
 \ $S \ = \ \bigcup_{n \geq 1} {\mathfrak S}_{P A^n}$ \  is a subgroup of 
$G_{k,1}$, and we have \ $G_{k,1} = S \cdot F_{k,1}$. 
Moreover, $S \cap F_{k,1} = \{ {\bf 1} \}$, hence we have a unique 
factorization.

The group $T_{k,1}$ has the subgroup \ 
$Z \ = \ \bigcup_{n \geq 1} {\mathbb Z}_{P A^n}$ \ and
we have \ $T_{k,1} = Z \cdot F_{k,1}$, with unique factorization for every
element of $T_{k,1}$.
\end{thm}
{\bf Proof.} 
Every element $\varphi$ of ${\mathfrak S}_{P A^n}$, as an element of
$G_{k,1}$, has finite domain and image codes that are the same: 
${\rm domC}(\varphi) = {\rm imC}(\varphi)$. Hence by Lemma \ref{lp_cap_F}, 
if $\varphi \in F_{k,1}$ then $\varphi = {\bf 1}$. Hence, 
$S \cap F_{k,1} = \{ {\bf 1} \}$, which implies uniqueness of the 
factorization as we saw in the beginning of the proof of Theorem 
\ref{factorization_lpV_F}.

To prove existence of the factorization we use the same factorization 
algorithm as in the proof of Theorem \ref{factorization_lpV_F}. Let 
$\varphi': P_1 \to Q_1$ represent any element of
$G_{k,1}$, where $P_1$ and $Q_1$ are finite maximal prefix codes. Then by
restriction we obtain a representation of the same element of the form
$\varphi: P_2 \to P A^n$, where it suffices to choose $n$ such that
$P A^n A^* \subseteq Q_1A^*$; since $P$ and $Q_1$ are finite maximal prefix
codes, such an $n$ exists. The remainder of the proof follows from the same
idea as for Theorem \ref{factorization_lpV_F}. We let $f: P_2 \to P A^n$ be
the (unique) element of $F_{k,1}$ determined by the finite maximal prefix
codes $P_2$ and $P A^n$, and let $\pi = \varphi \, f^{-1}$; then
${\rm domC}(\pi) = {\rm imC}(\pi) = P A^n$, hence 
$\pi \in {\mathfrak S}_{P A^n}$.

When $\varphi \in T_{k,1}$ the unique factorization $\varphi = \pi \, f$
satisfies $\pi = \varphi \, f^{-1} \in T_{k,1}$ (since 
$F_{k,1} \subset T_{k,1}$), hence $\pi \in T_{k,1} \cap S = Z$. 
 \ \ \ $\Box$

\medskip

Higman  (in \cite{Hig74}, Section 6) shows that the question whether a given
element of $G_{k,1}$ has finite order, is decidable.
The following theorem shows that every element of finite order of $G_{k,1}$
belongs to some subgroup ${\mathfrak S}_P$, for some finite maximal prefix
code $P$. Note that \ ${\rm domC}(\varphi) = {\rm imC}(\varphi) = P$ iff
$\varphi \in {\mathfrak S}_P$.

\begin{thm} \label{torsion_elements} \
Let $\Phi \in G_{k,1}$. Then $\Phi$ has finite order iff for some restriction
$\varphi$ of $\Phi$ we have
 \ ${\rm domC}(\varphi) = {\rm imC}(\varphi)$.
\end{thm}
{\bf Proof.}  \
If ${\rm domC}(\varphi) = {\rm imC}(\varphi) = P$ then
$\varphi \in {\mathfrak S}_P$, hence $\varphi$ has finite order.

Conversely, suppose that $\Phi$ is of finite order $r$, i.e.,
$\Phi^r(.) = {\rm id}(.)$ with $r > 0$, and $\Phi^i(.) \neq {\rm id}(.)$
for $0 \leq i < r$. By sufficiently restricting $\Phi$ we obtain maximal
finite prefix codes $P_0, P_1, \ldots, P_r \subset A^*$ such that for some
restriction $\varphi: A^* \to A^*$ of $\Phi$ we have \
 \ $P_0 \stackrel{\varphi}{\longrightarrow} P_1 $
 $\stackrel{\varphi}{\longrightarrow} \ \ldots \ $
$ \stackrel{\varphi}{\longrightarrow} \ P_{r-1} $
 $\stackrel{\varphi}{\longrightarrow} P_r$, and $\varphi(P_i) = P_{i+1}$ for
$i=0,1, \ldots, r-1$. Since $\varphi^r(.) = {\rm id}(.)$ it follows that
$P_0 = P_r$.

\smallskip

\noindent {\bf Claim:} For every $x \in P_0$, \
$C_x = \{ \varphi^i(x) : i = 0, 1, \ldots, r-1\}$ \ is a prefix code. \\
(Note: We only claim that no two $\varphi^i(x)$ are strict prefixes of each
other; we do not rule out that $\varphi^i(x) = \varphi^j(x)$ for some
$0 \leq i \neq j < r$.)

\smallskip

\noindent Proof of the Claim. If, by contradiction, we have
$\varphi^{\ell}(x) = x \, z$, for a non-empty word $z \in A^*$ and $x \in P_0$,
then for all $m \geq 0$ we have: \ $\varphi^{m \ell}(x) = x \, z^m$. This
implies that $\bigcup_{\ell=0}^{r-1} P_{\ell}$ contains words of arbitrarily
large length, which contradicts the fact that the prefix codes $P_{\ell}$ are
finite.

It follows that when $i>j$ then $\varphi^j(x) \in P_j$ cannot be a strict
prefix of $\varphi^i(x)$, since applying $\varphi^{-j}$ to
$\varphi^i(x) = \varphi^j(x) \, z$ yields $\varphi^{i-j}(x) = x \, z$.

Similarly, if we have $\varphi^i(x) = \varphi^j(x) \, z$ for a non-empty word
$z \in A^*$ and $x \in P_0$ and if $i<j$ then, applying $\varphi^{r-j}$ yields
$\varphi^{r+i-j}(x) = x \, z$, and the reasoning in the first paragraph (with
$\ell = r+i-j$) again yields a contradiction.
We conclude that $\varphi^i(x)$ and $\varphi^j(x)$ cannot be strict prefixes
of each other.  \ \ \ [End, Proof of Claim.]

\smallskip

\noindent The Claim implies that $\varphi(C_x) = C_x$ and that $C_x$ is a
cycle of $\varphi$. For each $x \in P_0$ we have a cycle $C_x$ as
above. For different $x \in P_0$ the corresponding cycles yield either
the same set or disjoint sets, i.e., for each $x, y \in P_0$, either
$C_x = C_y$ or $C_x \cap C_y = \emptyset$. So, $P_0$ is partitioned
into cycles of $\varphi$, hence $\varphi(P_0) = P_0$.
 \ \ \ $\Box$

\smallskip

As a consequence of Theorem \ref{torsion_elements} and Lemma
\ref{lp_cap_F} we recover a result of Brin and Squier \cite{BrinSqu}:

\begin{cor}
The group $F_{k,1}$ is torsion-free.
\end{cor}


\begin{thm} \label{PAm_equal_diff}
 \ If $P_1, P_2 \subset A^*$ are finite maximal prefix codes let
$S_i = \ \bigcup_{n \geq 0} {\mathfrak S}_{P_i A^n}$ \ for $i = 1$ or $2$.
We have: 

\smallskip

$S_1 = S_2$ \ iff \   
$\{ P_1A^n : n \geq 0\} \cap \{ P_2A^m : m \geq 0\} \neq \emptyset$,

\smallskip

\noindent 
When $S_1 \neq S_2$, the subgroup generated by $S_1 \cup S_2$ contains 
infinitely many elements of $F_{k,1}$.
\end{thm}
{\bf Proof.} \ 
If $P_1A^N = P_2A^M$ for some $M,N \geq 0$ then \ $S_1 = $
$\bigcup_{n \geq 0} {\mathfrak S}_{P_1 A^n} = $
$\bigcup_{n \geq 0} {\mathfrak S}_{P_1 A^N A^n}$, since 
${\mathfrak S}_{P_1 A^i} \hookrightarrow {\mathfrak S}_{P_1 A^N}$ when
$i \leq N$.  
Similarly, $\bigcup_{m \geq 0} {\mathfrak S}_{P_2A^M A^n} = S_2$.
Now, since $P_1A^N = P_2A^M$ we have \
$\bigcup_{n \geq 0} {\mathfrak S}_{P_1 A^N A^n} = $
$\bigcup_{m \geq 0} {\mathfrak S}_{P_2A^M A^n}$, hence $S_1 = S_2$.

\smallskip

In the other direction, under the condition
 \ $\{ P_1A^n : n \geq 0\} \cap \{ P_2A^m : m \geq 0\} = \emptyset$ \
we will prove that the subgroup of $G_{k,1}$ generated by $S_1$ and $S_2$
together contains some non-identity elements of $F_{k,1}$. Since $S_1$ and 
$S_2$ are torsion groups whereas $F_{k,1}$ is torsion-free, this implies that
$S_1 \neq S_2$.

\smallskip

\noindent
{\bf Claim.} \ There exist $n_0, m_0 \geq 0$ such
that \ $P_1A^{n_0} \cap P_2A^{m_0} \neq \emptyset$, \ and \   
$P_1A^{n_0} \neq P_2A^{m_0}$.
Moreover, there are $v_1 \in P_1A^{n_0} - P_2A^{m_0}$ and 
$v_2 \in P_2A^{m_0} - P_1A^{n_0}$ such that $v_1$ is a strict prefix of $v_2$.

\noindent Proof of the Claim: First, since each $P_1$ is a finite maximal
prefix code, every long enough word belongs to $P_1A^*$;  e.g.,
every word $w \in A^*$ of length $\geq {\rm max}\{ |p| : p \in P_1\}$
belongs to $P_1A^*$.  Therefore, for all $m$ large enough (e.g., all
$m \geq {\rm max}\{ |p| : p \in P_1\}$) we have $P_2A^m \subseteq P_1A^*$.
Let $m_0 \geq 0$ be such that $P_2A^{m_0} \subseteq P_1A^*$, and let us
consider the possible $n \geq 0$ such that
$P_2A^{m_0} \subseteq P_1A^nA^*$. For every $p_2u \in P_2A^{m_0}$ there exists
exactly one $p_1 v \in P_1A^n$ such that $p_1 v$ is a prefix of $p_2u$.
If $p_1 v \neq p_2u$, we can increase the length of $v$ (i.e., increase $n$)
to move $p_1 v$ closer to $p_2u$, until $p_1 v = p_2u$. Thus, there exists 
$n_0$ such that
 \ $P_1A^{n_0} \cap P_2A^{m_0} \neq \emptyset$.

Finally, $P_1A^{n_0} \neq P_2A^{m_0}$ by the hypothesis that
$\{ P_1A^n : n \geq 0\} \cap \{ P_2A^m : m \geq 0\} = \emptyset$.
Since $P_1A^{n_0} \neq P_2A^{m_0}$ and since $P_1A^{n_0}$ and $P_2A^{m_0}$ 
are finite maximal prefix codes, $P_1A^{n_0}$ and $P_2A^{m_0}$ are not strict 
subsets of each other. Hence there exist
$w_1 \in P_1A^{n_0} - P_2A^{m_0}$ and $v_2 \in P_2A^{m_0} - P_1A^{n_0}$.
Moreover, since $P_2A^{m_0} \subset P_1A^{n_0} A^*$ we have: For any 
$v_2 \in P_2A^{m_0}$ there exists $v_1 \in P_1A^{n_0}$ such that $v_1$ is a 
prefix of $v_2$. Since $v_2 \notin P_1A^{n_0}$, $v_1$ is a strict prefix of 
$v_2$.    \ \ \ [This proves the Claim.]

\medskip

We will now construct an element $\gamma_1 \in S_1$ whose $S_2 \cdot F_{k,1}$
factorization is of the form $\gamma_1 = \pi_2 \, f$ with $f \neq {\bf 1}$.
From this we obtain two elements $\gamma_1 \in S_1$ and $\pi_2 \in S_2$ such
that $\pi_2^{-1} \gamma_1 = f \in F_{k,1}$ with $f \neq {\bf 1}$.

Since $P_1A^{n_0} \cap P_2A^{m_0} \neq \emptyset$, there is 
$u_1 \in P_1A^{n_0} \cap P_2A^{m_0}$.
Using $u_1$ and the words $v_1$ and $v_2$ from the Claim, we now define \   
$\gamma_1 = (u_1|v_1)$; \ i.e., $\gamma_1$ is the permutation of $P_1A^{n_0}$
that transposes the two words $u_1$ and $v_1$, and fixes the rest of 
$P_1A^{n_0}$.

By Theorem \ref{PAm_factorization}, $\gamma_1 = \pi_2 \, f$ for a unique
$\pi_2 \in S_2$ and $f \in F_{k,1}$. The factorization algorithm given in the
proof of Theorem \ref{PAm_factorization} finds $\pi_2$ and $f$ by restricting 
$\gamma_1$ so that its image code becomes ${\rm imC}(\gamma_1) = P_2A^{m_0}$; 
the image codes of $f$ and of $\pi_2$ (not necessarily maximally extended), as 
well as the domain code of $\pi_2$, will also be $P_2A^{m_0}$.
The table of $\gamma_1$ is \ \ \

\smallskip

$ \gamma_1 \ = \ \left[ \begin{array}{cc l}
u_1 & v_1 \ & \ \ \ {\rm identity \ on}  \\
v_1 & u_1 \ & \ \ \ P_1A^{n_0} - \{u_1,v_1\}
\end{array}        \right]. $

\smallskip

\noindent By restricting so as to make ${\rm imC}(\gamma_1) = P_2A^{m_0}$
we obtain a table of the form 

\smallskip

$ \gamma_1 \ = \ \left[ \begin{array}{lll ll}
\ldots \ & \ u_1 z           \ & \ \ldots \ & \ v_1 \ & \ \ldots \\
\ldots \ & \ v_1 z \ (= v_2) \  & \ \ldots \ & \ u_1 \ & \ \ldots
\end{array}        \right].                  $

\smallskip

\noindent Here $z \in A^*$ is such that $v_2 = v_1 z$ and $z$ is non-empty
(recall that $v_1$ is a strict prefix of $v_2$).
Hence for this restriction of $\gamma_1, \pi_2$ and $f$ we have: \
${\rm domC}(\gamma_1) = {\rm domC}(f)$ contains $u_1 z$. But since
$u_1 \in P_2 A^{m_0}$  and since $P_2 A^{m_0}$ is a prefix code we find that
$u_1 z \notin P_2 A^{m_0} = {\rm imC}(f)$. Hence,
${\rm domC}(f) \neq {\rm imC}(f)$, therefore $f$ is not the identity. 
Since $F_{k,1}$ is torsion-free, the conclusion follows.
 \ \ \ $\Box$

\begin{thm} \label{PAm_conjugate}
 \ If $P_1, P_2 \subset A^*$ are finite maximal prefix codes let
$S_i = \ \bigcup_{n \geq 0} {\mathfrak S}_{P_i A^n}$ \ for $i = 1$ or $2$.
If \ $|P_1A^N| = |P_2A^M|$ \ for some
$N, M \geq 0$ then as subgroups of $G_{k,1}$, \
$S_1 = \theta^{-1} \ S_2 \ \theta$ \ for some $\theta \in F_{k,1}$.
\end{thm}
{\bf Proof.} \  Let $\theta$ be the element of $F_{k,1}$ such that 
$\theta: P_1A^N \to P_2A^M$; then $\theta$ can be restricted such that 
$\theta: P_1A^NA^n \to P_2A^MA^n$ for all $n \geq 0$.  
Then as subgroups of $G_{k,1}$, \ $S_1 = \theta^{-1} \ S_2 \ \theta$. 
 \ \ \ $\Box$

\bigskip

\noindent {\bf Element-specific factorizations}

\smallskip

\noindent For any element $\varphi \in G_{k,1}$ with ${\rm domC}(\varphi) = P$
and ${\rm imC}(\varphi) = Q$ we can apply Theorem \ref{PAm_factorization} to 
obtain the factorizations \ $\varphi(.) = \pi_Q \, f(.) = f \, \pi_P(.)$, 
where \ $f: P \to Q$ belongs to $F_{k,1}$, $\pi_Q \in {\mathfrak S}_Q$ and
$\pi_P \in {\mathfrak S}_P$. Moreover, \  
${\mathfrak S}_P = f^{-1} \, {\mathfrak S}_Q \, f$. Note that in this 
factorization, $\|f\|, \|\pi_Q\|, \|\pi_P\| \leq \|\varphi\|$.  

If $\varphi, \psi \in G_{k,1}$ are such that ${\rm domC}(\varphi) = P$, 
${\rm imC}(\varphi) = Q = {\rm domC}(\psi)$, and ${\rm imC}(\psi) = R$, 
then (since domain and ranges match)
we have the following multiplication formula for the factorization of 
$\psi \, \varphi(.)$. If $\varphi(.) = \pi_Q^{\varphi}\, f^{\varphi}$ and 
$\psi(.) = \pi_R^{\psi} \, f^{\psi}$ then  \   
$\psi \, \varphi(.) = \pi \, f(.)$, where  \  
$\pi \ = \ \pi_R^{\psi} \, f^{\psi} \pi_Q (f^{\psi})^{-1} \ \in $
${\mathfrak S}_R$, \ and \  
$f \ = \ f^{\psi} \, f^{\varphi} \ \in F_{k,1}$.

\bigskip

\noindent {\bf Questions left open:} \ What are all the torsion subgroups of 
$G_{k,1}$? What are all the torsion, non-torsion, or
torsion-free subgroups $S$ of $G_{k,1}$ for which there is a unique
factorization $G_{k,1} = S \cdot F_{k,1}$?
Are the groups \ $\bigcup_{n \geq 0} {\mathfrak S}_{P_1 A^n}$ \ and \ 
$\bigcup_{n \geq 0} {\mathfrak S}_{P_2 A^n}$ \ non-isomorphic if they do not
obey the conditions of Theorem \ref{PAm_conjugate}?


\section{Generators of {\it lp}$V$ and reversible computing}

\noindent
We are interested in the computation of elements of $V$ and of {\it lp}$V$ 
by circuits. For general information on circuits see \cite{Wegener, Handb}; 
good references on reversible circuits are 
\cite{FredToff, Toff80Memo, Toff80Conf, Shende}.
We will use the following fundamental results from the field of
{\it reversible computing}:

\smallskip

\noindent
$\bullet$ (V.~Shende, A.~Prasad, I.~Markov, J.~Hayes \cite{Shende}) \
Every {\em even} permutation of the set $\{0,1\}^n$ can be computed by a
circuit constructed only from bijective gates of type {\sc not, c-not, cc-not}.

The gates {\sc not, c-not, cc-not} are well known in the field of reversible
computing, and are defined as follows:

\smallskip

\noindent
{\sc not}: \ $x \in \{0,1\} \longmapsto \overline{x} \in \{0,1\}$ \ is the
usual negation operation;

\smallskip

\noindent {\sc c-not}: \ $(x,y) \in \{0,1\}^2 \longmapsto (x, y \oplus x) \in $
$\{0,1\}^2$ \ is the {\em controlled} {\sc not}, also called the
``Feynman gate'';

\smallskip

\noindent {\sc cc-not}: \ $(x, y, z) \in \{0,1\}^3 \longmapsto $
$(x, y, z \oplus (x \& y))$ $\in \{0,1\}^3$ \ is the {\em doubly controlled}
{\sc not}, with $\oplus$ denoting the usual exclusive {\sc or} (i.e.,
addition modulo 2), and $\&$ denoting the logical {\sc and} (i.e.,
multiplication modulo 2). The doubly controlled {\sc not} is usually called
the ``Toffoli gate'' \cite{FredToff, Toff80Memo, Toff80Conf}.

\smallskip

\noindent
$\bullet$ \cite{Shende} \ For any positive integer $n$ and any permutation
$\pi \in {\mathfrak S}_{2^n}$, the permutation \

\smallskip

$(x_1, \ldots, x_n, x_{n+1})$ $\in$ $\{0,1\}^{n+1}$
$ \ \longmapsto \ \pi(x_1, \ldots, x_n) \cdot x_{n+1}  \in \{0,1\}^{n+1}$ \

\smallskip

\noindent is {\it even}. Here, ``$\cdot$'' denotes concatenation. 
Indeed, one transposition $(u|v)$ (with $u,v \in \{0,1\}^n$) is now
replaced by $(u0|v0) \, (u1|v1)$ \ (i.e., two transpositions).

As a consequence, every odd permutation of $\{0,1\}^n$ can be computed
by a circuit that only makes use of bijective gates of type {\sc not, c-not,
cc-not}, and that uses an extra identity wire $x_{n+1} \mapsto x_{n+1}$.

\smallskip

\noindent
$\bullet$ (T.~Toffoli \cite{Toff80Memo, Toff80Conf}) An odd permutation of
$\{0,1\}^n$ cannot be computed by any circuit containing only bijective
gates with fewer than $n$ input-output wires. Hence for odd permutations,
the extra identity wire is necessary for bijective computing with a finite
collection of gate types.

\medskip

\noindent The above results have some interesting consequences for the group 
{\it lp}$V$:

\smallskip

First, the overall restriction operation for elements of the Thompson group 
$V$ (which replaces each $\varphi(x) = y$ by the pair $\varphi(x0) = y0$, 
$\varphi(x1) = y1$ for all $x$ in the domain of $\varphi$) now receives a 
very concrete interpretation for elements of {\it lp}$V$: For an element
$\varphi$ of {\it lp}$V$, the overall restriction is equivalent to adding an
{\it identity wire} at the ``bottom'' of the circuit (i.e., at the right-most
position for boolean variables).

\smallskip

Second, another consequence of the above concerns the generators of 
{\it lp}$V$.  Let $N, C, T$ be the partial maps $\{0,1\}^* \to \{0,1\}^*$ 
defined as follows, where $w \in \{0,1\}^*$ is any bitstring; \
$N: x_1w \mapsto \overline{x}_1w$, \
$C: x_1x_2w \mapsto x_1 \, (x_2 \oplus x_1) \, w$, \ and \
$T: x_1x_2x_3w \mapsto x_1x_2 \, (x_3 \oplus (x_2 \& x_1)) \, w$. These maps
are just the {\sc not, c-not, cc-not} gates, {\it applied only to the first
(left-most) bits} of a binary string. We leave $N, C, T$ undefined on bit 
strings that are too short.

\smallskip

Note that the engineering convention consists of using the same name
(e.g., ``not'', ``c-not'', etc.) for the same operation on different variables
in an sequence of variables. But this convention would not be correct in our
setting; e.g.,
negating the first bit in a string is different from negating the second bit.
In order to implement the operations $N, C, T$ on all bit-positions, i.e.,
in order to obtain the gate types {\sc not, c-not, cc-not} in the engineering
sense of the word, we also introduce the position transpositions
$\tau_{i,j}: \{0,1\}^* \to \{0,1\}^*$ (where $1 \leq i<j$), defined by

\smallskip

 \ \ \ \ \ $\tau_{i,j}: \ $
$u \, x_i \, v \, x_j \, w \ \longmapsto \ u \, x_j \, v \, x_i \, w$,

\smallskip

\noindent where $u, v, w \in \{0,1\}^*$, $|u| = i-1$, $|v| = j-1-i$; and we
leave $\tau_{i,j}(s)$ undefined when $|s| < j$.
Note that $\tau_{i,j}$ does not transpose a pair of words $\in \{0,1\}^*$,
but boolean variables (or positions within words).

\begin{pro} \label{nct_generators_lpV}
The group {\it lp}$V$ is generated by the set \
$\{N, C, T\} \, \cup \, \{\tau_{i,i+1} : 1 \leq i \}$.
More generally, {\it lp}$G_{k,1}$ is generated by \ $\Gamma_k \, \cup $
$\{\tau_{i,i+1} : 1 \leq i \}$ \ for some finite set $\Gamma_k$.

The finite alternating group ${\mathfrak A}_{2^n}$ (acting on $\{0,1\}^n$)
is generated by the set
 \ $\{N, C, T\} \cup \{\tau_{i,i+1} : 1 \leq i < n\}$.
\end{pro}
{\bf Proof.} The Proposition is immediate from the above observations, and
in particular the work \cite{Shende}.
 \ \ \ $\Box$

\bigskip

\noindent {\bf Application: An intuitive generating set for $V$}

\smallskip

The ${\it lp}V \cdot F$ factorization, together with the nice generating set
given above for  ${\it lp}V$ enables us to find a finite generating set for 
$V$ with a nice ``physical'' interpretation. It follows from the 
${\it lp}V \cdot F$ factorization that 
 \ $\{N, C, T\} \, \cup \, \{\tau_{i,i+1} : 1 \leq i \}$ 
$ \, \cup \ \{\sigma, \sigma_1\}$ is a generating set of $V$, where 
$\{\sigma, \sigma_1\}$ is the generating set of $F$ given in \cite{CFP}
with tables 

\smallskip

\hspace{1in} 
$\sigma \ = \ \left[ \begin{array}{lll}
00 & 01 & 1 \\
0 & 10 & 11
\end{array}    \right]  $, \ \ \ \ \  
$\sigma_1 \ = \ \left[ \begin{array}{llll}
0 & 100 & 101 & 11 \\
0 & 10 &  110 & 111
\end{array}    \right].   $

\smallskip

\noindent We will see that $\sigma_1$ is a ``controlled lowering'' of $\sigma$
(defined below). In \cite{BiThomps} we saw that $\sigma$ can be viewed 
as the ${\mathbb Z}$-shift $0^n 1 \mapsto 0^{n-1}1$, $01 \mapsto 10$, 
$1^n 0 \mapsto 1^{n+1} 0$, on the maximal prefix code $0^*01 \cup 1^*10$.  

Since $V$ is finitely generated, only a finite subset of
$\{\tau_{i,i+1} : 1 \leq i \}$ will be needed for generating $V$. 
Surprisingly, it turns out that in the presence of the other generators,
the Toffoli gate $T$ will not be needed for $V$. In detail we have:

\begin{pro} \label{generators_V}
The Thompson group $V$ is generated by the finite set \  
$\{N,\ C, \ \tau_{1,2}, \ \sigma, \ \sigma_1 \}$, \ where $N$ is the {\sc not} 
gate applied to the first wire, $C$ is {\sc c-not} (controlled {\sc not}, 
a.k.a.\ Feynman gate) applied to the first two wires,
$\tau_{1,2}$ is the transposition of the first two wires, and 
$\sigma, \ \sigma_1$ generate the Thompson group $F$.
\end{pro}
{\bf Proof.} We start out with the Higman generators $\kappa, \lambda, \mu, \nu$
of $V$ (see \cite{Hig74}), whose tables are   

\smallskip

$\kappa \ = \ \left[ \begin{array}{ll}
0 & 1 \\
1 & 0 
\end{array}   \right]  $, \ \  
$\lambda \ = \ \left[ \begin{array}{lll}
00 & 01 & 1 \\
00 & 1 & 01
\end{array} \right] $, \ \  
$\mu \ = \ \left[ \begin{array}{lll}
0 & 10 & 11 \\
10 & 0 & 11
\end{array} \right] $, \ \ 
$\nu \ = \ \left[ \begin{array}{llll}
00 & 01 & 10 & 11 \\
00 & 10 & 01 & 11
\end{array} \right] $. 

\smallskip

\noindent We see that $\kappa = N$ and $\nu = \tau_{1,2}$. For 
$\lambda$ and $\mu$ we apply the ${\it lp}V \cdot F$ factorization 
algorithm, which leads to

\smallskip

$\lambda(.) \ = \ \left[ \begin{array}{llll}
00 & 01 & 10 & 11 \\
00 & 10 & 11 & 01
\end{array} \right] \ \cdot \    
\left[ \begin{array}{llll}
00 & 010 & 011 & 1  \\
00 & 01 & 10 & 11 
\end{array} \right](.) $

\smallskip

\noindent The right factor belongs to $F$, hence it is generated by 
$\{ \sigma, \sigma_1\}$. It is easy to check that the first factor is
equal to \ $\tau_{1,2} \cdot C(.)$; recall that $C$ 
has the table \ 
$\left[ \begin{array}{llll}
00 & 01 & 10 & 11 \\
00 & 01 & 11 & 10 
\end{array} \right] \ = \ 
\left[ \begin{array}{lll}
0 & 10 & 11 \\
0 & 11 & 10 
\end{array} \right] $.
  
\smallskip

\noindent
A similar calculation leads to a factorization \ $\mu(.) \ = \ $
$\tau_{1,2} \cdot C \cdot \tau_{1,2} \cdot N \cdot \tau_{1,2} \cdot f(.)$,
for some $f \in F$.                     \ \ \ $\Box$

\bigskip

\noindent {\bf The lowering operation} 

\smallskip

The following operation, inspired by circuits, gives further insight into 
${\it lp}V$ and $F$. For any integer $d>0$ we define

\smallskip

$\varphi \in G_{k,1} \ \longmapsto \ (\varphi)_d \in G_{k,1}$ \ by

\smallskip

$(\varphi)_d(zx) = z \ \varphi(x)$ \ for all $z \in A^d$, \
$x \in {\rm Dom}(\varphi)$.

\smallskip

\noindent Recall that $A^d$ is the set of all words of length $d$ over $A$.
It is easy to see that for each $d>0$, the operation
$\varphi \to (\varphi)_d$ is an {\em endomorphism} of $G_{k,1}$, which is
injective but not surjective;
it is also an endomorphism of ${\it lp}G_{k,1}$, of $F_{k,1}$, of $T_{k,1}$,
and of ${\it lp}T_{k,1}$.

\smallskip

The circuit interpretation of the operation $\varphi \to (\varphi)_d$ is
that the ``gate'' $\varphi$ is lowered by $d$ positions in the circuit  
through 
the introduction of $d$ identity wires on top of the ``gate'' $\varphi$ 
(i.e., at the left end of the list of input variables). While $\varphi$ is 
applied to the boolean variables $x_1, x_2, \ldots$, the lowered gate will be
applied to the variables $x_{d+1}, x_{d+2}, \ldots$.  This is commonly done in 
circuits, as it allows the designer to place gates at any place in the 
circuit. In electrical engineering, traditionally no distinction is 
made between a gate; e.g., the {\sc c-not} operation and its lowerings are all
just called ``{\sc c-not} gates''.
The lowering operation is an important link between circuits and their
representation by groups or monoids of functions. 

\smallskip

The lowering operation can be expressed in terms of the transpositions
$\tau_{i,i+1}$, although the formula depends on the length $\ell$ of the 
longest word $\in A^*$ appearing in the table of $\varphi$.
We have: \ \ $(\varphi)_d(.) \ = \ \pi^{-1} \ \varphi \ \pi(.)$, \ where 
$\pi$ is the following permutation of bit positions:

\smallskip

\noindent If if $d+1 > \ell$ then 
$ \ \ \pi(.) \ = \ $
$\left( \begin{array}{cccc cccc}
\hspace{-.1in} 1 & 2 & \ldots & \ell & d+1 & d+2 & \ldots & d+\ell \\
\hspace{-.1in} d+1 & d+2 & \ldots & d+\ell & 1 & 2 & \ldots & \ell
\end{array} \hspace{-.09in} \right)$.

\smallskip

\noindent If $d+1 \leq \ell$ then 
$ \ \ \pi(.) \ = \ $
$\left( \begin{array}{cccc cccc}
\hspace{-.1in} 1  & 2  & \ldots & d & d+1 & d+2 & \ldots & d+\ell \\
\hspace{-.1in} 1+\ell & 2+\ell & \ldots & d+\ell & 1 & 2 & \ldots & \ell
\end{array} \hspace{-.09in} \right)$.

\medskip

When we write elements of $G_{k,1}$ as elements of the Cuntz algebra
${\mathcal O}_k$ (according to \cite{BiThomps} and \cite{Nekrash}), we 
see that the lowering operation is an endomorphism of ${\mathcal O}_k$ 
given by the formula

\smallskip

$\gamma \in {\mathcal O}_k \ \longmapsto \ (\gamma)_d = $
$\sum_{z \in A^d} z \, \gamma \, {\overline z} \ \in \ {\mathcal O}_k$.

\smallskip

Note that all transpositions of variables (or wires) $\tau_{i,i+1}$ are
obtained from the transposition of variables $\tau_{1,2}$ by \
$\tau_{i,i+1} = (\tau_{1,2})_{i-1}$. So, the lowering operations, together
with a finite set of elements of ${\it lp}G_{k,1}$, yields a generating
set of ${\it lp}G_{k,1}$.
For $G_{k,1}$ we already saw that the transpositions of variables are 
redundant as generators (since $G_{k,1}$ is finitely generated), but that 
the use of the transpositions of variables shortens the word length; 
we will see (Theorem \ref{VwordsCircuits}) that when the transpositions 
of variables are added to a finite
generating set of $G_{k,1}$ then the word length becomes approximately the
same as the bijective-circuit complexity.
For $F_{k,1}$ it will be interesting to consider generating sets of the
form $\Gamma \cup \ \bigcup_{d \geq 1} (\Gamma)_d$, where $\Gamma$ is any
finite generating set of $F_{k,1}$, and where
$(\Gamma)_d = \{ (f)_d : f \in \Gamma \}$ (see the open problems at the
end of Section 6). 

\smallskip

We can define the {\bf controlled lowering operation}; we fix any string 
$c \in A^*$, called the ``control string'' and define

\smallskip

$\varphi \in G_{k,1} \ \longmapsto \ (\varphi)_c \in G_{k,1}$ \ by

\smallskip

$(\varphi)_c(cx) \ = \ c \ \varphi(x)$ \ for all $x \in {\rm Dom}(\varphi)$, 
and

\smallskip

$(\varphi)_c(p \alpha) \ = \ p \alpha$ \ where $p <_{{\rm pref}} c$, \  
$\alpha \in A$, \ $p \alpha \not\leq_{{\rm pref}} c$;

\smallskip

\noindent so here $p$ is any strict prefix of
$c$ (i.e., $p \neq c$), and $\alpha \in A$ is such that $p \alpha$ is
{\em not} a prefix of $c$. 
So, domC$((\varphi)_c)  \ = \ c \cdot {\rm domC}(\varphi) \ \cup \ $
$\{p \alpha : p <_{{\rm pref}} c, \ \alpha \in A, \ $
$p \alpha \not\leq_{{\rm pref}} c\}$, \ and \ imC$((\varphi)_c)  \ = \ $
$c \cdot {\rm imC}(\varphi) \ \cup \ $
$\{p \alpha : p <_{{\rm pref}} c, \ \alpha \in A, \ $
$p \alpha \not\leq_{{\rm pref}} c\}$.

It is easy to see that for each $c \in A^*$, the operation
$\varphi \to (\varphi)_c$ is an {\em endomorphism} of $G_{k,1}$, which is
injective but not surjective;
it is also an endomorphism of ${\it lp}G_{k,1}$, of $F_{k,1}$, of $T_{k,1}$,
and of ${\it lp}T_{k,1}$.
In Cuntz algebra notation, the operation takes the form

\smallskip

$\gamma \in {\mathcal O}_k \ \longmapsto \ (\gamma)_c \ = \ $
$c  \, \gamma \, {\overline c}  \ + \ $
$\sum_{p, \alpha } p \alpha \, {\overline {p \alpha}} \ \in \ {\mathcal O}_k$.

\smallskip

\noindent where $p$ ranges over the strict prefixes of $c$ and $\alpha$ 
ranges over the letters of $A$ such that $p \alpha$ is not a prefix of $c$.
In ${\mathcal O}_k$ the controlled lowering operation is a multiplicative
endomorphism, but it is not additive.

Observe that for the generators $\{ \sigma, \sigma_1 \}$ of $F$ seen before,
the second generator is the controlled lowering of the first with control
string 1; this explains our notation for $\sigma_1$.


\section{Generalized word problem, distortion of $F_{k,1}$ and 
  ${\it lp}G_{k,1}$ in $G_{k,1}$ }

\begin{pro} 
Over any finite generating set of $G_{k,1}$ the generalized word problem 
of $F_{k,1}$ in $G_{k,1}$ can be decided in {\em cubic} deterministic time. 

Similarly, over any finite generating set of $G_{k,1}$ the generalized word 
problem of ${\it lp}G_{k,1}$ in $G_{k,1}$, and more generally, the 
generalized word problem of \ 
$\bigcup_m {\mathfrak S}_{PA^m}$ \ in $G_{k,1}$ (for any finite maximal 
prefix code $P$) can be decided in {\em cubic} deterministic time.
\end{pro}
{\bf Proof.} By Proposition 4.2 of \cite{BiThomps}, if $\varphi$ is given 
by a word of length $n$ over a finite generating set of $G_{k,1}$ then
a table for $\varphi$ (not necessarily maximally extended) can be computed 
in time $O(n^3)$. By Proposition 3.5 of \cite{BiThomps}, the length $n$ 
provides a linear upper bound on the size of this table. Also, every table 
entry has length $\leq c \ n$. 
More precisely, the table has the form $((x_1,y_1), \dots, (x_N,y_N))$, 
where $|x_i|, |y_i|, N \leq c \ n$ (for some constant $c \geq 1$). The sets 
$\{x_1, \dots, x_N\}$ and $\{y_1, \dots, y_N\}$ are maximal prefix codes,
and $\varphi(x_i) = y_i$ for $i = 1, \dots, N$. 

To check whether $\varphi$ belongs to $F_{k,1}$ we first sort the table 
according to the input entries, with respect to dictionary order; more 
precisely, we sort the pairs of the table $((x_1,y_1), \dots, (x_N,y_N))$ 
according the $x$-coordinates, in time $\leq O(n^2 \, {\rm log} \, n)$; 
indeed, there are $O(n \, {\rm log} \, n)$ sorting steps, and since each 
word has length $\leq c \ n$, each word comparison takes time $O(n)$. Then 
we check whether the resulting $x$-sorted table is now also in sorted form 
regarding the $y$-coordinates; this takes quadratic time, as there are 
$O(n)$ words of length $O(n)$.    

To check whether $\varphi$ belongs to ${\it lp}G_{k,1}$ we check, in time
$\leq O(n^2)$, that $|x_i| = |y_i|$ for $i = 1, \ldots, N$. And to check
whether $\varphi$ belongs to \ $\bigcup_m {\mathfrak S}_{PA^m}$ we note 
first that $P$ is a fixed finite maximal prefix code, independent of 
$\varphi$. We restrict $\varphi$ so that every table entry receives length
$\geq {\rm max}\{|p| : p \in P\}$. This multiplies the table size of $\varphi$
by a constant, at most (since $P$ and ${\rm max}\{|p| : p \in P\}$ are fixed). 
So we can assume that each $x_i$ and $y_i$ in the table of $\varphi$ has a 
prefix in $P$.
Now for $x_1$, find the prefix $p_1 \in P$ of $x_1$, so $x_1 = p_1s_1$ for 
some $s_1 \in A^*$, and let $m_0 = |s_1|$. Thus, $\varphi$ belongs to \ 
$\bigcup_m {\mathfrak S}_{PA^m}$ iff $\varphi \in {\mathfrak S}_{PA^{m_0}}$. 
So, we now write each $x_i$ and each $y_i$ in the form 
$p \, s$ with $p \in P$ and check 
that $|s| = m_0$; this holds (for all $s$ obtained) iff $\varphi \in$
${\mathfrak S}_{PA^{m_0}}$.  Checking this takes time $\leq O(n^2)$.
 \ \ \ $\Box$

\medskip

Since ${\it lp}G_{k,1} \cap F_{k,1} = \{ {\bf 1} \}$, we have the following
equivalence: \ $w = {\bf 1}$ (as elements of $G_{k,1}$) iff 
$w \in {\it lp}G_{k,1}$ and $w \in F_{k,1}$. 
Thus, the word problem of $G_{k,1}$ reduces (by a one-to-one linear-time 
reduction) to the conjunction of the generalized word problem
of ${\it lp}G_{k,1}$ in $G_{k,1}$ and the generalized word problem of $F_{k,1}$
in $G_{k,1}$.  (Here, the reduction function is just the identity map.)
The same is true with ${\it lp}G_{k,1}$ replaced by 
$\bigcup_m {\mathfrak S}_{PA^m}$ (for any chosen finite maximal prefix 
code $P$).

Hence, the deterministic (or nondeterministic, or co-nondeterministic) time
complexity of the word problem of $G_{k,1}$ is a lower bound for the 
deterministic (respectively nondeterministic, or co-nondeterministic) time 
complexity of the generalized word problem of 
$\bigcup_n {\mathfrak S}_{PA^m}$ in $G_{k,1}$ or the generalized word 
problem of $F_{k,1}$ in $G_{k,1}$, or both. 
More formally, we have the following:

\begin{defn} 
We say that a language (or decision problem) $L$ is {\em as hard as} 
coNP iff there is a coNP-complete problem $L_0$ such that for every function 
$t(.)$ that is a deterministic time complexity lower bound for infinitely 
many instances of $L_0$ we have: Some function $\geq c \cdot t(.))$ is a 
deterministic time complexity lower bound for infinitely many instances 
of $L$ (for some constant $c > 0$).
\end{defn} 

\begin{defn} \label{defn_wordlength}
Let $G$ be a group with generating set $A$. Suppose every generator
$\alpha \in A$ has been assigned a ``length'' $|\alpha| \in {\mathbb N}$.
Typically, if $A$ is finite then $|\alpha| =1$ for all $\alpha \in A$.
For the position transpositions $\tau_{i,j}$ ($1 \leq i < j$) we take 
$|\tau_{i,j}| = j$.

The length of a word $w = a_1 \ldots a_n$ over $A$ is defined by \
$|w| = \sum_{j=1}^n |a_j|$.

The {\em word length} $|g|_A$ of $g \in G$ over $A$ is defined to be the
shortest length of any word (over $A$) that represents $g$.
\end{defn}
For a group with generating set $A$ we often say ``a word over $A$'' when
we actually mean ``a word over $A \cup A^{-1}$''; we will also use the
notation $A^{\pm 1}$ for $A \cup A^{-1}$.

\begin{pro} \label{as_hard_as_coNP}
Let $\Gamma_{k,1}$ be a finite generating set of $G_{k,1}$ but suppose that
elements of $G_{k,1}$ are given over the infinite generating set \  
$\Gamma_{k,1} \cup \{\tau_{i,j}: 1 \leq i<j \}$. Then the generalized word 
problem, in $G_{k,1}$, of either $F_{k,1}$ or ${\it lp}G_{k,1}$, or both, 
is as hard as coNP.

Similarly, if $P \in A^*$ is a finite maximal prefix code then the generalized 
word problem, in $G_{k,1}$, of either $F_{k,1}$ or \, 
$\bigcup_n {\mathfrak S}_{PA^m}$, or both, is as hard as coNP. 
Moreover, these problems are in coNP.
\end{pro}
{\bf Proof.} The word problem of $G_{k,1}$ over the generating set
$\Gamma_{k,1} \cup \{\tau_{i,j}: 1 \leq i<j \}$ is coNP-complete
\cite{BiCoNP}. The hardness then follows from the above conjunctive 
reduction.  
 \ \ \ $\Box$

\begin{defn} \label{defn_distortion}
Let $G_1$ be a group with generating set $A_1$, and let $G_2$ be a subgroup
of $G_1$ with generating set $A_2$. A function 
$f: {\mathbb N} \to {\mathbb N}$ is called a {\bf distortion} function for 
$G_2$ within  $G_1$, with respect to the generators $A_1$, respectively $A_2$,
iff for all $g_2 \in G_2$: \ $|g_2|_{A_2} \leq f(|g_2|_{A_1})$. 

{\em The} distortion function of $G_2$ within  $G_1$, with respect to the
generators $A_1$, respectively $A_2$, is the smallest distortion function.
\end{defn} 

\begin{pro}
If we use {\em finite} generating sets for both $G_{k,1}$ and $F_{k,1}$ 
then $F_{k,1}$ has linear distortion in $G_{k,1}$.
\end{pro}
{\bf Proof.} For any element $g \in G_{k,1}$ we have \   
$\|g\| \leq c_1 \, |g|_G$, by Proposition 3.5 of \cite{BiThomps}; here, 
$\|g\|$ is the table size of the element $g$, $c_1$ is a positive constant, 
and $|g|_G$ is the word length of $g$ over some chosen, fixed finite 
generating set of $G_{k,1}$. 
By Theorem 2.5 of \cite{CFP}, \ $|g|_F \leq c_2 \, \|g\|$, where $c_2$ is a
positive constant, and $|g|_F$ is the word length of $g$ over some chosen,
fixed finite generating set of $F_{k,1}$. Hence, $|g|_F \leq c_1c_2 \, |g|_G$,
so the distortion of $F_{k,1}$ in $G_{k,1}$ is linear. 
 \ \ \ $\Box$

\bigskip

\noindent {\bf Problems left open:} \ 

\noindent 1. \ \ Over the generating set
$\Gamma_{k,1} \cup \{\tau_{i,j}: 1 \leq i<j \}$ of $G_{k,1}$, are the
generalized word problems of the subgroups $F_{k,1}$, ${\it lp}G_{k,1}$, 
and $\bigcup_n {\mathfrak S}_{PA^n}$ each coNP-complete?

\smallskip

\noindent 2. \ \ We saw that ${\it lp}V$ is generated by \  
$\{ N, C, T\} \ \cup \ \{\tau_{i,i+1}: 1 \leq i \}$.
What is the {\em distortion} of ${\it lp}V$ (over this generating set)
within the Thompson group $V$ (with $V$ over the generating set 
$\Gamma_V \ \cup \ \{\tau_{i,i+1}: 1 \leq i \}$, where $\Gamma_V$ is 
any finite generating set of $V$)?  We will see in the next Section that 
this distortion has a close connection to the relation between different
kinds of bijective circuits.

\noindent 3. \ \ We saw that $F$ is generated by a two-element set 
$\{ \sigma, \sigma_1\}$, and hence also by 
$\{ (\sigma)_d, (\sigma_1)_d : d > 0\}$. What is the {\em distortion} of 
$F$ within $V$, when $F$ is taken over the generating set 
$\{ (\sigma)_d, (\sigma_1)_d : d > 0\}$, and $V$ is taken over the 
generating set $\Gamma_V \, \cup \{\tau_{i,i+1}: 1 \leq i \}$?


\section{Complexity of $F$ and of the factorization of $V$}

We saw that in the factorization $\varphi = \pi \, f$ with 
$\pi \in {\it lp}V$ and $f \in F$, the table sizes of $\pi$ and 
of $f$ can be exponentially larger than the table size of $\varphi$. 
We will now investigate the circuit complexity of $\pi$ and $f$, compared to
that of $\varphi$.
We will also show that if $f \in F$ then the circuit complexity of $f^{-1}$ 
is not much higher than the circuit complexity of $f$; in other words, the 
elements of the Thompson group $F$ do not have much computational asymmetry
(and in particular, they cannot be one-way functions).   
And we will show that some problems in $V$ are coNP-complete or 
$\#{\mathcal P}$-complete; in particular, the problem of finding the  
${\it lp}V \cdot F$ factorization is $\#{\mathcal P}$-complete.
In this section we focus on the Thompson groups $V$ and $F$, but the results
could easily be extended to $G_{k,1}$ and $F_{k,1}$. 


\subsection{Circuit complexity and Thompson groups}

Since an element $\varphi \in V$ is a partial function mapping bitstrings to
bitstrings, it is natural to view $\varphi$ as a boolean function, to be 
computed by a boolean circuit. However, unless $\varphi \in {\it lp}V$, the 
inputs and the outputs of $\varphi$ do not have a fixed length. 
So the traditional concept of a combinational boolean circuit cannot be 
applied directly to elements of $V$.

Let $\varphi: P \to Q$ be a bijection between finite maximal prefix codes 
$P,Q \subset \{0,1\}^*$, representing an element of $V$.
We will use ternary logic over the alphabet $\{0,1,\bot\}$, where $\bot$ is
a new letter used for padding bitstrings.
Let $m$ is the length of the longest bitstring in $P$, and let $n$ is the
length of the longest bitstring in $Q$. 
We define $\varphi^{\bot}: \{0,1,\bot\}^m \to \{0,1,\bot\}^n$ as follows:

\smallskip

  For $p \in P$, \ \
  $\varphi^{\bot}(p \, \bot^{m -|p|}) \ = \ $
    $ \varphi(p) \ \bot^{n - |\varphi(p)|}$.

\smallskip

  For $x \in \{0,1,\bot\}^m - \{ p \, \bot^{m -|p|} : p \in P\}$
  we let \ $\varphi^{\bot}(x) \ = \ \bot^n$.

\smallskip

\noindent We will use the notation \ 

\smallskip

$P^{\bot} \ = \ \{ p \, \bot^{m -|p|} : p \in P\}$
$ \ = \ P\bot^* \ \cap \ \{0,1,\bot\}^m$, \ where 
$m = {\rm max}\{ |p| : p \in P\}$; 

\smallskip

$Q^{\bot} \ = \ \{ q \, \bot^{n -|q|} : q \in Q\}$ 
$ \ = \ Q\bot^* \ \cap \ \{0,1,\bot\}^n$, \ where 
$n = {\rm max}\{ |q| : q \in Q\}$.

\smallskip

\noindent
We call $P^{\bot}$, $Q^{\bot}$, and $\varphi^{\bot}$ the {\bf padding} of 
$P$, $Q$, respectively $\varphi$.
Note: For $\varphi \in V = G_{2,1}$, the padding $\varphi^{\bot}$ 
is {\em not} to be viewed as an element of {\it lp}$G_{3,1}$.

We observe that for the restrictions to $P^{\bot}$ or to $Q^{\bot}$ we have \ 
$(\varphi^{\bot}|_{P^{\bot}})^{-1} = (\varphi^{-1})^{\bot}|_{Q^{\bot}}$; 
the restriction $\varphi^{\bot}|_{P^{\bot}}$ is bijective (but 
$\varphi^{\bot}$ is not bijective in general).
When ${\rm imC}(\varphi) = {\rm domC}(\psi)$ we also have \
$(\psi \varphi)^{\bot} = \psi^{\bot}\varphi^{\bot}$.

\medskip

To compute the function 
$\varphi^{\bot}: \{0,1,\bot\}^m \to \{0,1,\bot\}^n$ we consider
ternary-logic combinational circuits with gates over the alphabet
$\{0,1,\bot\}$. We assume that a finite, computationally universal set of
ternary logic gates has been chosen; we ignore the details since they only
affect the circuit complexity by a constant multiple. We also use the
(unbounded) set of wire-swap operations $\tau_{i,i+1}$. 

For such a circuit, the {\bf size of the circuit} is defined to be the 
number of gates together with the number of wires (links between gates
or between gates and inputs or outputs). Note that a ``lowered gate'' 
$(\gamma)_d$ (i.e., the gate $\gamma$ applied to the wires $d+1, d+2$, etc.,
as defined at the end of Section 5) is counted as one gate (independently of 
$d$). Also, in a circuit each wire-crossing $\tau_{i,i+1}$ will be counted as 
one gate (independently of $i$). Note that here we are talking about circuit
size, not about word length.

\medskip

\noindent {\bf Remarks}: \\
(1) \ The idea of padding with $\bot$ works for the Thompson-Higman group
$G_{k,1}$ in general, by using $(A \cup \{\bot\})$-valued logic.
The gates that we use include the wire-crossings $\tau_{i,i+1}$. \\   
(2) \ In \cite{BiDistor} we will follow another, more algebraic, approach
for defining circuit complexity of elements of $V$. We embed
$V$ into a certain finitely generated partial transformation {\it monoid} $M$
acting on $\{0,1,\bot\}^*$, and we take the word-length of $\varphi$ in $M$
as the circuit complexity of $\varphi$. We will prove in \cite{BiDistor} that
there are monoids $M$ that, over certain generators, can ``simulate'' logic 
gates, and that in such monoids word-length is closely related to circuit 
complexity.

\medskip

Since the functions  $\varphi: P \to Q$ considered here are elements of $V$,
hence bijective, it is natural to also introduce {\bf bijective 
$\{0,1,\bot\}$-valued circuits}.
A $\{0,1,\bot\}$-valued circuit is said to be bijective iff the gates that 
make up the circuit are the wire-swap operations $\tau_{i,i+1}$ ($i \geq 1$),  
and a set of gates derived from the elements of some fixed finite generating 
set $\Gamma_V$ of $V$. The latter means, more precisely, that the gates derived
from $\Gamma_V$ are of the form $((\gamma)_d)^{\bot}$ where each $\gamma$ is a 
restriction of an element of $\Gamma_V$. Recall that $(\gamma)_d$ (for 
$d \geq 0$) is the lowering of $\gamma$ (defined at the end of Section 5).

\smallskip

In this paper, unless we specifically mention ``bijective'' or
``$\{0,1,\bot\}$-valued'', the word ``circuit'' will refer to a general
boolean circuit (not necessarily bijective).

\medskip
 
\noindent {\it Comparison between $\{0,1,\bot\}$-valued circuits and boolean
(i.e., $\{0,1\}$-valued) circuits:}

In the general (not necessarily bijective) case, a $\{0,1,\bot\}$-valued
circuit can be simulated by a traditional binary-logic circuit (e.g., 
by encoding the ternary values $0, 1, \bot$ by the binary strings 
$00, 11, 01$ respectively, with $10$ also serving as a code for $\bot$). 
Thus, there is no essential difference between general 
$\{0,1,\bot\}$-valued circuits and general boolean circuits.

However, bijective $\{0,1,\bot\}$-valued circuits have greater generality 
than bijective boolean circuits. First, bijective boolean circuits have 
input-output functions belonging to  ${\it lp}V$ only; on the other hand,
input-output functions of bijective $\{0,1,\bot\}$-valued circuits are the
paddings of all the elements of $V$. Also, the input-output function of a 
bijective $\{0,1,\bot\}$-valued circuit is only bijective as a function 
$P^{\bot} \to Q^{\bot}$, not as a function  
$\{0,1,\bot\}^m \to \{0,1,\bot\}^n$, whereas the input-output function of
a bijective boolean circuit is a permutation of $\{0,1\}^m$ for some $m$. 

Moreover, even for $\varphi \in {\it lp}V$ the smallest size of a bijective 
$\{0,1,\bot\}$-valued circuit computing $\varphi$ is the word length of 
$\varphi$ over the generators of $V$ (as we shall show in Theorem
\ref{VwordsCircuits} below), whereas the smallest size of a bijective 
boolean circuit computing $\varphi$ is the word length of $\varphi$ over 
the generators of ${\it lp}V$. 
Thus, we will see that the relation between the two bijective circuit sizes 
is approximately the {\it distortion} of ${\it lp}V$ within $V$. Here the 
generating set of $V$ is $\Gamma_V \cup \{ \tau_{i,i+1} : 1 \leq i\}$ for 
any finite generating set $\Gamma_V$ of $V$, and ${\it lp}V$ is generated by 
$\{N,C,T\} \cup \{ \tau_{i,i+1} : 1 \leq i\}$ as seen before. 
Finding the distortion of ${\it lp}V$ within $V$ is one of our open problems
mentioned at the end of Section 6. Theorem \ref{VwordsCircuits} below will
give a precise connection between the two kinds of bijective circuit 
sizes, word lengths in $V$ or in ${\it lp}V$, and the distortion of 
${\it lp}V$ in $V$.  

\begin{defn}
Let $\Gamma_{k,1}$ be a finite generating set of $G_{k,1}$, and
let $w$ be a word over the generating set $\Gamma_{k,1} \, \cup $
$ \{\tau_{i,i+1} : 1 \leq i \}$ of $G_{k,1}$. The {\bf length} of 
$w = a_1 \ldots a_n$ is \, $|w| \ = \ \sum_{j=1}^n |a_j|$, where
$|a_j| = 1$ if $a_j \in \Gamma_{k,1}$ and $|a_j| = i+1$ if  
$a_j = \tau_{i,i+1}$.   

For $\varphi \in G_{k,1}$, the {\bf word length} of  $\varphi$ over
$\Gamma_{k,1} \, \cup \,  \{\tau_{i,i+1} : 1 \leq i \}$ is the shortest 
length of any word (over the above generators) that represents $\varphi$.
\end{defn}

Observe that $\tau_{i,i+1}$ is counted differently for circuit size than 
for word length ($\tau_{i,i+1}$ is counted as 1 in circuit size but as
$i+1$ in the word length). 

The following definition compares bijective padded circuits for elements of
${\it lp}V$ with boolean bijective circuits. Note that in this 
definition we only consider circuits for computing elements of ${\it lp}V$.  

\begin{defn}
An {\em unpadding cost function} from bijective $\{0,1,\bot\}$-valued 
circuits to bijective binary circuits is any function 
$U: {\mathbb N} \to {\mathbb N}$ such that the following holds: For all 
$\varphi \in {\it lp}V$ and any bijective $\{0,1,\bot\}$-valued circuit 
of size $m$ for $\varphi$ there exists a bijective binary circuit of size
$\leq U(m)$ for $\varphi$.

{\em The} unpadding cost function $u(.)$ is the minimum unpadding cost 
function. 
\end{defn}

\begin{defn}
Two functions $f_1, f_2: {\mathbb N} \to {\mathbb N}$ are said to be {\bf 
linearly related} iff there are constants $c_0, c_1, c_2$, all $\geq 1$, such 
that for all $n \geq c_0$: \    
$f_1(n) \leq c_1 \, f_2(c_1 n)$ \ and \ $f_2(n) \leq c_2 \, f_1(c_2 n)$.

The functions $f_1, f_2$ are said to be {\bf polynomially related} iff there
are constants $c_0, c_1, c_2$, all $\geq 1$, such that for all 
$n \geq c_0$: \ $f_1(n) \leq c_1 \, f_2(c_1 n^{c_1})^{c_1}$ \ and \     
$f_2(n) \leq c_2 \, f_1(c_2 n^{c_2})^{c_2}$.  
\end{defn}

The following theorem motivates $\{0,1,\bot\}$-valued circuits, as well as 
the concept of word length over the generating set $\Gamma_V \, \cup $ 
$ \{\tau_{i,i+1} : 1 \leq i \}$ for $V$. It again motivates our use of
the infinite set $ \{\tau_{i,i+1} : 1 \leq i \}$ for generating $V$,
inspite of the fact that $V$ is finitely generated. It also reinforces 
the connection between Thompson groups and bijective (``reversible'') 
computing, seen before. 
In \cite{BiDistor} we will generalize Theorem \ref{VwordsCircuits} to a 
connection between general circuit size and the word size in the 
``Thompson monoids'' (the latter being a 
generalization of the Thompson-Higman groups to monoids).
We will only state the Theorem for $V$ and for binary or $\{0,1,\bot\}$-valued
circuit, although it could easily be generalized to $G_{k,1}$.
First a lemma:

\begin{lem} \label{length_in_table} \    
Let $a_n, \ldots, a_1 \in V$ be given by table and let $\ell$ be the length 
of the longest words in the tables of $a_n, \ldots, a_1$.
Then $a_n, \ldots, a_1$ have restrictions $\alpha_n, \ldots, \alpha_1$,
respectively, such that 

$\bullet$ \ \ ${\rm domC}(\alpha_{j+1}) = {\rm imC}(\alpha_j)$ for
$n >j \geq 1$, and 

$\bullet$ \ \ all words in the tables of $\alpha_j$ ($n \geq j \geq 1$) 
have lengths \ $\leq n \, \ell$.
\end{lem}
{\bf Proof.} 
For $\varphi \in V$ we will describe the table of $\varphi$ as a set of
input-output pairs, of the form $\{(u_i, v_i) : i = 1, \ldots, I\}$.
We also use tables to represent elements of $V$ in non-maximally extended
form; we will mention explicitly when we assume maximal extension.

\medskip

\noindent {\bf Claim.}
{\it \ Let $a_i \in V$ (for $i = 1, \ldots, n$) be given by tables
$\{(x_j^{(i)},y_j^{(i)}) : j = 1, \ldots, r\}$, not necessarily in reduced
form. Thus
 \ ${\rm domC}(a_i) = \{ x_j^{(i)} : j = 1, \ldots, r\}$, and
 \ ${\rm imC}(a_i) = \{ y_j^{(i)} : j = 1, \ldots, r\}$.
We assume that ${\rm domC}(a_{i+1}) = {\rm imC}(a_i)$ for
$i = 1, \ldots, n-1$
Let $\ell$ be an upper bound on the length of all the words in
 \ $\bigcup_{i=1}^n {\rm domC}(a_i) \cup {\rm imC}(a_i)$.
For $j = 1, \ldots, r$, let $S_j = \{ s_{j,k} : 1 \leq k \leq |S_j|\}$ be
a finite maximal prefix code over $\{0,1\}$.

Then ${\alpha}_i$ (for $i = 1, \ldots, n$), defined by the table
$\{(x_j^{(i)}s_{j,k}, \ y_j^{(i)}s_{j,k}) : 1 \leq j \leq r, $
  $1 \leq k \leq |S_j|\}$, is a restriction of $a_i$ satisfying:

$\bullet$ \ \ ${\rm domC}(\alpha_{i+1}) = {\rm imC}(\alpha_i)$ \ for
$n > i \geq 1$, and

$\bullet$ \ \ all words in the tables of $\alpha_i$ ($n \geq i \geq 1$)
have lengths \ $\leq $
${\rm max}\{|s_i| : 1 \leq i \leq r\} + \ell$.
}

\smallskip

\noindent Proof of the Claim: 
 Since \ $\{ x_j^{(i+1)} : j = 1, \ldots, r\} = $
${\rm domC}(a_{i+1}) = {\rm imC}(a_i) = \{ y_j^{(i)} : j = 1, \ldots, r\}$,
it follows immediately that \
$\{ x_j^{(i+1)} s_{j,k} : 1 \leq j \leq r, 1 \leq k \leq |S_j|\} = $
$\{ y_j^{(i)} s_{j,k} : 1 \leq j \leq r, 1 \leq k \leq |S_j|\}$.
Hence ${\rm domC}(\alpha_{i+1}) = {\rm imC}(\alpha_i)$.

Also, $|x_j^{(i)} s_{j,k}| \leq \ell + {\rm max}\{|s_i| : 1 \leq i \leq r\}$
(and similarly for $|y_j^{(i)} s_{j,k}|$), hence we have the claimed length
bound. This proves the Claim.

\medskip

Let us now prove Lemma \ref{length_in_table} by induction on $n$. The Lemma
is obvious when $n=1$. Given $a_i \in V$ (for $i = n, \ldots, 1$, with
$n \geq 2$), we use the Lemma by induction for $a_{n-1}, \ldots, a_1$.
So we can assume that ${\rm domC}(a_{i+1}) = {\rm imC}(a_i)$ for
$n-1 > i \geq 1$, and all words in the tables of $a_i$ ($n-1 \geq i \geq 1$)
have lengths $\leq (n-1) \, \ell$. Let us denote
${\rm domC}(a_{i+1}) = {\rm imC}(a_i)$ by $P_i$; so we have
 \ $P_0 \stackrel{a_1}{\to} P_1 \stackrel{a_2}{\to} \ \ldots \ $
$ \stackrel{a_{n-2}}{\longrightarrow} P_{n-2}$
$ \stackrel{a_{n-1}}{\longrightarrow} P_{n-1}$ .

We consider the product $a_n \cdot (a_{n-1} \ldots a_1)$.
We will find a restriction $\alpha_n$ of $a_n$, and a restriction \
$\overline{a_{n-1} \ldots a_1}$ \ of $a_{n-1} \ldots a_1$, such that
${\rm domC}(\alpha_n) = {\rm imC}(\overline{a_{n-1} \ldots a_1})$.
We also want to restrict $a_{n-1}, \ldots, a_1$ to functions
$\alpha_{n-1}, \ldots, \alpha_1$ such that
${\rm domC}(\alpha_{i+1}) = {\rm imC}(\alpha_i)$ for $n-1 > i \geq 1$.
Two cases arise:

\smallskip

\noindent Case 1: Every word in $P_{n-1}$ has length $\geq \ell$.

By the assumptions of Lemma \ref{length_in_table}, every word in the table
of $a_n$ has length $\leq \ell$. Therefore, all we need to do to obtain
$\alpha_n$ and $\alpha_{n-1}, \ldots, \alpha_1$ is to restrict $a_n$ so
that ${\rm domC}(\alpha_n)$ becomes $P_{n-1}$.
No restriction of $a_{n-1} \ldots a_1$ is needed, i.e., $\alpha_i = a_i$ for
$i = n-1, \ldots, 1$, and \
$\overline{a_{n-1} \ldots a_1} = a_{n-1} \ldots a_1$. Hence the longest word
in $P_0, P_1, \ldots, P_{n-1}$ has length $\leq (n-1) \, \ell$.

The longest word in the table of $\alpha_n$ has length
$\leq {\rm max}\{|p| : p \in P_{n-1}\} + \ell$, by the
Claim (applied to $\alpha_n$ and $a_{n-1} \ldots a_1$). Hence the longest
word in the table of $\alpha_n$ has length
$ \leq (n-1) \, \ell + \ell = n \, \ell$.

\smallskip

\noindent Case 2: Some word in $P_{n-1}$ has length $< \ell$.

We restrict $a_{n-1}$ so as to make all words in ${\rm imC}(\alpha_{n-1})$
have length $\geq \ell$, as follows. For any $y_j \in P_{n-1}$ with
$|y_j| < \ell$ we consider the finite maximal prefix code
$S_j = \{0,1\}^{\ell - |y_j|}$. We restrict $a_{n-1}$ in such a way that
$y_j$ is replaced by $y_j \cdot S_j$, i.e., $P_{n-1}$ becomes \
$(P_{n-1} - \{y_j\}) \, \cup \, y_j \cdot S_j$. Note that all words in
$y_j \cdot S_j$ have length $\ell$.  After every word in $P_{n-1}$
of length $< \ell$ has been replaced, let $\overline{P}_{n-1}$ be the
resulting finite maximal prefix code. Now we apply the Claim in order to
restrict all of $a_{n-1}, \ldots, a_1$. As a result, each $\alpha_i$ (for
$i = n-1, \ldots, 1$) receives a table with words of length \
$\leq (n-1) \, \ell + {\rm max}\{|s| : s \in \bigcup_j S_j\}$
$\leq (n-1) \, \ell + \ell = n \, \ell$. Note that in these restrictions, the
length of the words in $P_{n-1}$ only increases for the very short words
(namely, words of length $< \ell$ are replaced by words of length $\ell$).
Hence, after restriction, the words in $\overline{P}_{n-1}$ still
have length $\leq (n-1) \, \ell$.

Next we restrict $a_n$, as in case 1. Since after restriction, the words in
$\overline{P}_{n-1}$ have length $\leq (n-1) \, \ell$, the longest word in
the table of $\alpha_n$ has length $ \leq (n-1) \, \ell + \ell = n \, \ell$.
 \ \ \ $\Box$

\begin{thm} \label{VwordsCircuits} \ 
{\bf (1)} For the elements $\varphi \in V$, the minimum {\em size} of 
bijective $\{0,1,\bot\}$-valued circuits that compute $\varphi$ is 
polynomially related to the {\em word length} of $\varphi$ in $V$ (over 
$\Gamma_V \, \cup $ $ \{\tau_{i,i+1} : 1 \leq i \}$). 
More precisely, there are constants $c_1, c_2 >0$ such that \     
$s_{\varphi} \leq c_1 \, |\varphi|_{_V}^2$, and  \   
$|\varphi|_{_V} \leq c_2 \, s_{\varphi}$, where $|\varphi|_{_V}$ is word 
length of $\varphi \in V$ over $\Gamma_V \, \cup $
$ \{\tau_{i,i+1} : 1 \leq i \}$, and $s_{\varphi}$ is the 
$\{0,1,\bot\}$-valued bijective circuit size of $\varphi$.      

\smallskip

\noindent {\bf (2)} For the elements $\varphi \in {\it lp}V$, the minimum 
size of bijective binary circuits is polynomially related to the word length 
of $\varphi$ in ${\it lp}V$ (over 
$\{ N,C,T\} \, \cup \, \{\tau_{i,i+1} : 1 \leq i \}$). 
More precisely, the word length $|\varphi|_{_{{\it lp}V}}$ over 
$\{ N,C,T\} \, \cup \, \{\tau_{i,i+1} : 1 \leq i \}$, and the binary 
circuit size $b_{\varphi}$ of $\varphi$ satisfy:
 \ \ $b_{\varphi} \leq c_1 \, |\varphi|_{_{{\it lp}V}}^2$, and 
 \ $|\varphi|_{_{{\it lp}V}} \leq c_2 \, b_{\varphi}$ (for some constants
$c_1, c_2 >0$). 

\smallskip

\noindent {\bf (3)} The {\bf distortion} function $d(.)$ of ${\it lp}V$ 
in $V$ (over the generators mentioned above), and the {\bf unpadding cost}
function $u(.)$ for bijective $\{0,1,\bot\}$-valued circuits, are 
polynomially related.  More precisely, for some constants $c, c'>0$ and for 
all $x > 0$ we have 
 \ $u(x) \leq c \ d(c \, x)^2$ \ and \ 
$d(x) \leq c' \ u(c' \, x^2)$.
\end{thm}
{\bf Proof.} 
(1) For $\varphi \in V$ let $s_{\varphi}$ be the circuit size of $\varphi$
over $\Gamma_V^{\pm 1} \, \cup \, \{\tau_{i,i+1}: 1 \leq i \}$, let 
$|\varphi|$ be the word length over 
$\Gamma_V^{\pm 1} \, \cup \, \{\tau_{i,i+1}: 1 \leq i \}$.   

\smallskip

\noindent $\bullet$ Proof that \ $s_{\varphi} \leq c_1 \, |\varphi|^2$ 
(for some constant $c_1>0$): \    
Let $w = a_1 \ldots a_n$ be a shortest word that represents $\varphi$, where 
$a_j \in \Gamma_V^{\pm 1} \, \cup \, \{\tau_{i,i+1}: 1 \leq i \}$ 
for $1 \leq j \leq n = |w| = |\varphi|$. 
We restrict the generators $a_j$ as in Lemma \ref{length_in_table} so that 
they can be composed, and $a_j$ will only have bitstrings of length 
$\leq c_1 \, |\varphi|$ in its table. Thus, the word $w$ becomes a 
$\{0,1,\bot\}$-valued circuit consisting of the $n = |\varphi|$ operations 
$a_j$ ($1 \leq j \leq n$), and each $a_j$ has $\leq c_1 \, |\varphi|$ wires; 
so the circuit for $\varphi$ has size $\leq c_1 \, |\varphi| \, |\varphi|$. 

\smallskip

\noindent $\bullet$ Proof that \ $|\varphi| \leq c_2 \, s_{\varphi}$
(for some constant $c_2 >0$): \ 
Consider any smallest bijective circuit $C$ over 
$\Gamma_V^{\pm 1} \, \cup \, \{\tau_{i,i+1}: 1 \leq i \}$, of size
$s_{\varphi}$, computing $\varphi$.
This circuit is a sequence $(a_1, \ldots, a_n)$ where $n = |C| = s_{\varphi}$. 
Each $a_j$ is either of the form $\tau_{i,i+1}$, or $a_j$ is the padding of a 
restriction of $(\gamma)_d$ with $d \geq 0$ and 
$\gamma \in \Gamma_V^{\pm 1}$. Hence, the sequence $(a_1, \ldots, a_n)$ is a 
word of length $s_{\varphi}$ representing $\varphi$.
To obtain a word over $\Gamma_V^{\pm 1} \, \cup \, \{\tau_{i,i+1}: 1 \leq i \}$
we express the lowering operation in terms of position transpositions, as at 
the end of Section 5; then $(\gamma)_d$ becomes 
$\pi^{-1} \, \gamma \, \pi$, where 
$\pi$ is the composition of $\leq 2 \, m$ position transpositions of the 
form $\tau_{i,d+i}$ or $\tau_{i,m+i}$. Here $m$ is the length of 
the longest word in any of the tables for the elements $\gamma \in \Gamma_V$; 
since $\Gamma_V$ is fixed and finite, $m$ is a constant.
A transposition $\tau_{i,d+i}$ can be written as the composition of
$< 2d$ transpositions of the form $\tau_{j,j+1}$. 
Let $((\gamma_j)_{d_j} : j = 1, \ldots, J)$ be the list of all the lowered 
gates that occur in the circuit $C$; then \ $\sum_{j=1}^J d_j < |C|$, since
for each $(\gamma_j)_{d_j}$ there are $d_j$ wires in $C$ that are output wires
of other gates (or that are inputs of $C$), and that are counted as 
part of the size of $C$.   
Thus we obtain a word of length $< c_2 \, s_{\varphi}$ (for some constant 
$c_2 > 0$), representing $\varphi$. 

\smallskip

\noindent (2) The proof is very similar to the proof of (1).

\smallskip

\noindent (3) By (1) and (2) and by the definition of distortion we have: \   
$c \ \sqrt{b_{\varphi}} \ \leq \ |\varphi|_{_{{\it lp}V}} \ \leq \ $
$d(|\varphi|_{_V}) \ \leq \ d(c' \ s_{\varphi})$, hence \ 
$b_{\varphi} \ \leq \ c'' \ d(c' \ s_{\varphi})^2$. Here, $c, c', c''$ are 
constants. By the definition of the unpadding cost function it follows that 
 \ $u(x) \leq c'' \ d(c' \ x)^2$.  

Also by (1) and (2) and by the definition of the unpadding cost function we 
have: \ $c \, |\varphi|_{_{{\it lp}V}} \ \leq \ b_{\varphi} \ \leq \ $
$u(s_{\varphi}) \ \leq \ u(c' \ |\varphi|_{_V}^2)$, hence \ 
$|\varphi|_{_{{\it lp}V}} \ \leq \ c'' \ u(c' \ |\varphi|_{_V}^2)$. Here, 
$c, c', c''$ are constants. By the definition of distortion it follows 
that \ $d(x) \leq c'' \ u(c' \, x^2)$.
 \ \ \ $\Box$

\begin{thm} \label{inverting_f} 
Consider any element of the Thompson group $F$, represented by a bijection
$f: P \to Q$ that preserves the dictionary order, where $P$ and $Q$ are 
finite maximal prefix codes.
If $f$ can be computed by a $\{0,1,\bot\}$-valued circuit of size $s$ then
$f^{-1}: Q \to P$ can be computed by a combinational circuit size
 \ $m(m+1) \, s + O(m^2 n)$, where $m = {\rm max}\{|p| : p \in P\}$ and
$n = {\rm max}\{|q| : q \in Q\}$.

Moreover, a circuit for $f^{-1}$ can be found from a circuit for $f$ 
deterministically in polynomial time in terms of $s, m, n$.
\end{thm}
{\bf Proof.} Suppose $f^{\bot}(x \bot^{m-|x|}) = y \bot^{n-|y|}$, and 
$y \bot^{n-|y|}$ is given. The idea for inverting $f^{\bot}$ is simple: 
Since $f$ preserves the dictionary order we can find $x \bot^{m-|x|}$ 
by adapting the classical {\em binary search} algorithm. This algorithm is 
usually used for searching in a sorted array; but it works in a similar way 
for inverting any order-preserving map. 

A few technical details have to be discussed before we give an algorithm 
for inverting $f^{\bot}$.
For many strings $z \in \{0,1,\bot\}^n$ there is no inverse image under 
$f^{\bot}$; in that case our inversion algorithm will output $\bot^m$. For 
example, $z \not \in \{0,1\}^*\bot^*$ has no inverse. Also, since the range 
of $f^{\bot}$ is $Q^{\bot} \cup \{\bot^n\}$ where $Q$ is a finite maximal 
prefix code we have: 
If $y \bot^{n-|y|}$ has an inverse then there is no inverse for any strict 
prefix of $y$; i.e., if $w \in \{0,1\}^*$ is a strict prefix of $y$ then 
$w \bot^{n-|w|}$ has no inverse. On the other hand, for every 
$v \in \{0,1\}^n$, there exists exactly one prefix $y$ of $v$ such that 
$y \bot^{n-|y|}$ has an inverse. Similarly, for every $u \in \{0,1\}^m$,
there exists exactly one prefix $x$ of $u$ such that 
$f^{\bot}(x \bot^{m-|x|})$ $\neq \bot^n$.      

The binary search can be pictured on the complete binary tree with vertex
set $\{0,1\}^{\leq m}$,  with root $\varepsilon$ (the empty word), and 
leaves $\{0,1\}^m$; the children of a vertex $v \in \{0,1\}^{< m}$ are
$v0$ and $v1$. The search uses a variable vertex $v$, initialized to 
$\varepsilon$, and proceeds from $v$ to $v0$ or $v1$, until success, or 
until $v$ becomes a leaf.  

\medskip

\noindent {\tt Algorithm} \ (for inverting $f^{\bot}$)

\smallskip

\noindent For any input $z \in \{0,1,\bot\}^n$ it is easy to check whether 
$z \not\in \{0,1\}^*\bot^*$; in that case the output is $\bot^m$.
Assume from now on that the input is of the form $y \bot^{n-|y|}$, with 
$y \in \{0,1\}^*$, $|y| \leq n$. 
Let $v \in \{0,1\}^{\leq m}$ be a bitstring, initialized to $v = \varepsilon$ 
(the empty string).
Below, $<_{_d}$ and $>_{_d}$ refer to the dictionary order.

\smallskip

\noindent {\tt repeat the following until the exit command} \\ 
{\tt begin}  

find the prefix $x$ of \ $v \, 1 \, 0^{m-1-|v|}$ \ (if $|v| < m$) or the 
prefix $x$ of $v$ (if $|v| = m$), such that  

 \ \ \ \ \ $f^{\bot}(x \bot^{m-|x|}) \neq \bot^n$
 \ \ (this is done by trying all prefixes of $v \, 1 \, 0^{m-1-|v|}$,  
      respectively of $v$);  

{\tt if} \ $f^{\bot}(x \bot^{m-|x|}) = y \bot^{n-|y|}$ \ {\tt then} 
 return \ $x \bot^{m-|x|}$ \ as output, and {\tt exit};

{\tt if} \ $|v| =m$ \ {\tt then} return $\bot^m$ as output, and {\tt exit};

{\tt if} \ $f^{\bot}(x \bot^{m-|x|}) <_{_d} \ y \bot^{n-|y|}$  \ and $|v|<m$ 
  {\tt then} replace $v$ by $v1$; 

{\tt else} (i.e., when $f^{\bot}(x \bot^{m-|x|}) >_{_d} \ y \bot^{n-|y|}$ and
 $|v|<m$) \ replace $v$ by $v0$;

\noindent {\tt end (of repeat loop)};

\smallskip

\noindent {\tt end (of algorithm).}

\medskip

Let us show that this program can be implemented by an acyclic circuit of 
size \   $m(m+1) s + O(m^2 n)$. 
The loop of the program is executed at most $m$ times, and each 
execution of the loop will be implemented as one of $m$ stages of the
complete circuit. 

Each execution of the loop takes at most $(m+1)s + O(mn)$ gates. 

The first part of the loop (namely, to ``find the 
prefix $x$'') requires $m$ copies of the circuit of $f^{\bot}$, each of which 
is followed by $O(n)$ gates to check equality with $\bot^n$, followed by 
a tree of $n$ {\sc and}-gates.
Moreover, recall that when we want to apply the same type of gate to different
variables (wires) we need to permute wires (using bit position transpositions 
$\tau_{i,j}$). Similarly, permutations may need to be applied to the output
wires of a gate. This adds at most a constant number of operations for each
gate.
So the first part of the loop uses $ms+ O(mn)$ gates. 

The first {\tt if} condition requires another copy of the circuit of 
$f^{\bot}$, followed by $O(n)$ gates to compare the result with 
$y \bot^{n-|y|}$ for equality and to check for $<_{_d}$ or $>_{_d}$ in the 
dictionary order. The $<_{_d}$- or $>_{_d}$-comparison of two strings of the 
same length can be done by a finite automaton, reading both strings in 
parallel from left to right; if the inputs are restricted to strings of 
length $n$, this automaton can then be turned into a {\it prefix circuit} 
(of Ladner and Fischer \cite{LadnerFischer}). The Ladner-Fischer
circuit consists of $\leq 4n$ copies of a gate that implements the (fixed) 
transition function of the finite automaton. The prefix circuit uses fan-out 
$< n$; however, there is also a bounded-fan-out design for the prefix 
circuit, using just $< 9n$ gates (see p.\ 205 of \cite{LakshDhall}). 
Moreover, applying the same gate to different variables first requires some
permutations of wires; this introduces a constant factor (since gates have 
a fixed number of input-output wires, so only a fixed number of wires are 
permuted back and forth).
Checking whether $|v| = m$ is equivalent to checking absence of $\bot$, which
requires $O(n)$ gates.
So overall the if-statements require $s + O(n)$ gates. 

\smallskip

Finally, the above description amounts to a polynomial-time procedure for 
producing the circuit that implements $(f^{\bot})^{-1}$.
 \ \ \ $\Box$

\medskip

As a consequence, $(f^{\bot})^{-1}$ is not much harder to compute than 
$f^{\bot}$ itself. So, without need to define the concept of a one-way 
function in detail we can conclude that for any reasonable definition of 
``one-way function'' we have:

\begin{cor}
The Thompson group $F$ does not contain any one-way functions.
\end{cor}


Recall that in our algorithm for finding the ${\it lp}V \cdot F$ factorization 
of $\varphi \in V$, the element $\varphi$ is first restricted so as to make
${\rm imC}(\varphi) = \{0,1\}^n$.  We will show next that this restriction
does not increase circuit complexity much, and that we can find a circuit 
for certain restrictions. On the other hand, we saw in Theorem 
\ref{coNP_circuits}(3) that the opposite operation, namely finding the 
maximal extension, is hard. 

\begin{lem} \label{makingImC_n}
Every $\varphi \in V$ has a restriction $\Phi$ such that 
${\rm imC}(\Phi) = \{0,1\}^n$, and such that the circuit size of $\Phi$ is 
only polynomially larger than the circuit size of $\varphi$.

More precisely, assume $\varphi^{\bot}$ has a circuit of size $s$, with
$m$ input variables and $n$ output variables (over $\{0,1,\bot\}$). Then the 
restriction $\Phi$ with ${\rm imC}(\Phi) = \{0,1\}^n$ has a circuit of 
size $\leq 4s \, (n+m+1)$, with $n+m$ input variables and $n$
output variables. Moreover, such a circuit for $\Phi^{\bot}$ can be found from
the given circuit for $\varphi^{\bot}$ deterministically in polynomial time
(as a function of $s, m, n$); i.e., there is a polynomial-time reduction from
the problem of finding a circuit for $\Phi$ to the problem of finding a 
circuit for $\varphi$ (for $\{0,1,\bot\}$-valued circuits).
\end{lem}
{\bf Proof.} We can view $\varphi$ as a bijection $P \to Q$ where
$P,Q \subset \{0,1\}^*$ are finite maximal prefix codes. Since the circuit
for $\varphi^{\bot}$ has $n$ output variables, we have
$n = {\rm max}\{|y| : y \in Q\}$.
Let $\Phi: P_1 \to \{0,1\}^n$ be the restriction of $\varphi$ with image code
$\{0,1\}^n$, where $P_1$ is the finite maximal prefix code obtained when 
$\varphi$ is restricted to make the image code $\{0,1\}^n$. Let 
$m = {\rm max}\{|x| : x \in P\}$. Thus, all words in the
finite maximal prefix code $P_1$ have length $\leq n+m$. We now construct a
$\{0,1, \bot\}$-valued circuit for $\Phi^{\bot}$, with $n+m$ input variables 
and $n$ output variables. On an input $x \in \{0,1, \bot\}^{n+m}$ the circuit 
behaves as follows:

\smallskip

\noindent
$\bullet$ If $x \notin \{0,1\}^*\bot^*$, the output is $\Phi(x) = \bot^n$. 

To check whether $x \notin \{0,1\}^*\bot^*$ we consider all $n+m-1$ pairs 
$(x_i,x_{i+1})$ of neighboring input variables 
(for $i=1, \ldots, n+m-1$) and check whether any of them have values 
$(x_i,x_{i+1}) = (\bot,0)$ or $= (\bot,1)$, using $n+m-1$ gates. 
To produce the output $\bot^n$ in that case, the $n+m-1$ gates above feed into 
a tree of {\sc or} gates whose output is 1 iff $(\bot,0)$ or $(\bot,1)$ occurs 
anywhere in input pairs.
The {\sc or}-tree and the output $\bot^n$ require $< 2(n+m)$ gates. 
Thus so far we have $< 3(n+m)$ gates in total. 

\smallskip

\noindent
$\bullet$ If $x \in \{0,1\}^*\bot^*$, since $\Phi^{\bot}$ has $m+n$ input wires,
we write $x = u \, \bot^i$ with $u \in \{0,1\}^{n+m-i}$. We look at each prefix 
$p$ of $u = pz$, in order of increasing length $|p| = 0, 1, \ldots, m+n-i$, 
and feed $p \, \bot^{m-|p|}$ into $\varphi^{\bot}$.  

\smallskip

\noindent 
- If $\varphi^{\bot}(p \, \bot^{m-|p|}) = \bot^n$, we ignore $p$ and look 
at the next prefix of $u$.

\smallskip

\noindent 
- If $\varphi^{\bot}(p \, \bot^{m-|p|}) = q \, \bot^{m-|q|}$ \ for some 
$q \in \{0,1\}^*$, we conclude that
 $\varphi(p) = q$ and $\varphi(u) = \varphi(pz) = qz$. 
Hence, if $|z| = n - |q|$ we produce the output 
$\Phi(u \, \bot^i) = qz \in \{0,1\}^n$; so $\Phi$ agrees with $\varphi$ and has
imC$(\Phi) = \{0,1\}^n$.
If $|z| \neq n - |q|$ we produce the output $\Phi(u \, \bot^i) = \bot^n$.
(No new prefixes of $u$ will be considered.)  

\smallskip

In the above construction, the circuit of $\varphi^{\bot}$ is repeated $m+n+1$
times, since an input of length $m+n$ has $\leq m+n+1$ prefixes. So this part 
of the circuit has size $s \, (m+n+1)$.
We need another $3n \, (m+n+1)$ gates to combine the outputs of the $(m+n+1)$
copies of the $\varphi^{\bot}$-circuit: If one of the $\varphi^{\bot}$-circuits 
produces an output in $\{0,1\}^n$, that output has to be the final output; if 
all the copies of the $\varphi^{\bot}$-circuit produce $\bot^n$, then $\bot^n$ 
should be the final output.

Finally, the total circuit for $\Phi^{\bot}$ has size 
 \ $\leq 3(n+m) + s \, (m+n+1) + 3n \, (m+n+1)$ $ \leq (s+3n+3)(m+n+1)$
 $ \leq 4s \, (m+n+1)$ \ (the last ``$\leq$'' holds since 
$s \geq m+n$ and $m,n \geq 1$). 
The above description of the construction of a circuit for $\Phi^{\bot}$ is a
deterministic algorithm whose running time is a polynomial in $s$.
 \ \ \ $\Box$

\bigskip

An immediate consequence if Lemma \ref{makingImC_n} is the following:

\begin{cor} \label{f_imC_01m}
Assume $f \in F$ has a $\{0,1,\bot\}$-valued circuit of size $\leq s$, with 
$m$ input variables and $n$ output variables. Then the restriction of $f$  
with ${\rm imC}(f) = \{0,1\}^n$, i.e., the restriction of $f$ that makes $f$
a rank function ${\rm rank}_P(.)$, has a $\{0,1,\bot\}$-valued circuit of 
size $\leq 4s \, (m+n+1)$.  
\end{cor}
In other words, representing elements of $F$ by rank functions does not
lead to a large increase in circuit complexity.


\begin{pro} \label{complexity_of_factors}
Let $\varphi$ be an element of the Thompson group $V$, and let
$\varphi = \pi \cdot f$ be its ${\it lp}V \cdot F$ factorization.
Let $\varphi: P \to \{0,1\}^n$ be a representation of $\varphi$ by a 
bijection from a finite maximal prefix code $P \subseteq \{0,1\}^{\leq m}$ 
onto $\{0,1\}^n$.  Suppose that the rank function of $P$ can be computed 
by a circuit of size $\leq s$.

Then $f$ has circuit complexity $\leq s$, and the circuit complexities of
$\varphi$ and of $\pi$ differ by at most \ $m(m+1) \, s + O(m^2 n)$.  
The circuit complexities of $\varphi^{-1}$ and $\pi^{-1}$  differ by
at most  $s$.  

Moreover, the circuits for $f$ and $\pi$ can be found in deterministic
polynomial time.          
\end{pro}
{\bf Proof.} We apply our ${\it lp}V \cdot F$ factorization algorithm. Since
$\varphi$ already has imC$(\varphi) = \{0,1\}^n$, we have 
$f = {\rm rank}_P(.)$, and hence by assumption, $f$ has circuit complexity
$\leq s$. 
 
To obtain a circuit for $\pi = \varphi \, f^{-1}$ we use Theorem
\ref{inverting_f} to obtain a circuit for $f^{-1}$ of size  
$\leq m(m+1) \, s + O(m^2 n)$, where $m = {\rm max}\{|p| : p \in P\}$;
then we compose the circuit for $f^{-1}$ with the circuit for $\varphi$.

To obtain a circuit for $\varphi = \pi \cdot f$ from a circuit for $\pi$
we just compose the circuit for $f$ and the circuit for $\pi$. 
 \ \ \ $\Box$

\medskip

A consequence of Proposition \ref{complexity_of_factors} is the following. 
If an element of $V$ has a representation $\varphi: P \to \{0,1\}^n$ (for some
$n>0$), and if $P$ is a finite maximal prefix code with easy rank function,
then $\varphi$ and of $\pi$ have similar circuit complexities; $\varphi^{-1}$
and $\pi^{-1}$ also have similar circuit complexities. Thus we have:

\begin{cor}
If there exists a one-way bijection $\varphi: P \to \{0,1\}^n$ (for
some $n>0$), where $P$ is a finite maximal prefix code with easy rank function,
then there exists a one-way permutation $\pi$ of $\{0,1\}^n$.
\end{cor}


\subsection{coNP-complete and $\#{\mathcal P}$-complete problems in the
Thompson groups}

The following coNP-completeness results are similar to the well-known
coNP-completeness of questions about circuits, except that here we deal
with circuits that compute bijections, in the sense defined at the beginning
of this section.

\begin{thm} \label{coNP_circuits}
The following decision problems are {\bf coNP-complete}: \\
(1) Given two $\{0,1,\bot\}$-valued bijective circuits, do they compute the
 same element of $V$?  \\
(2) Given two $\{0,1,\bot\}$-valued bijective circuits for computing elements
 $\psi, \varphi \in V$, is $\psi$ the maximal extension of $\varphi$?  \\
(3) Given a $\{0,1,\bot\}$-valued bijective circuit, does it compute the
 identity element of $V$?     

The problems remain coNP-complete when the given 
$\{0,1,\bot\}$-valued circuits are general (not necessarily bijective).
\end{thm}
{\bf Proof.} Let us first check that these problems are in coNP. Problems
(1), (3) and (4) are variants of the classical circuit equivalence problem.
For problem (2), we can check in coNP whether $\psi$ and $\varphi$ represent
the same element of $V$. To check in NP whether $\psi$ is {\em not} maximally
extended, guess entries $(x0,y0), (x1,y1)$ in the table of $\psi$; the lengths
of $x$ and $y$ are no larger than the size of the given circuit for $\psi$,
and the fact that $(x0,y0)$ and $(x1,y1)$ are in the table of $\psi$ can be
checked rapidly using the circuit for $\psi$.

Hardness: Problem (3) is a special case of (2) and of (1) (letting $\psi$ be
the identity map with domain and image codes consisting of just the empty
word), so (2) and (1) are at least as hard as (3).
The hardness of (3) is a consequence of the fact that the word problem
of $G_{3,1}$ over the generating set $\Gamma_{3,1} \cup $
$\{\tau_{i,j} : 1 \leq i < j\}$ is coNP-complete (proved in \cite{BiCoNP}),
and the fact (proved in Theorem \ref{VwordsCircuits}) that every word over
$\Gamma_{3,1} \cup $ $\{\tau_{i,j} : 1 \leq i < j\}$ has a
$\{0,1,\bot\}$-valued
circuit whose size is linearly bounded by the size of the word.
 \ \ \ $\Box$

\bigskip

Proposition \ref{complexity_of_factors} shows that under certain conditions 
the ${\it lp}V \cdot F$ factorization is easy to find. The next Theorems show 
that in general, finding the ${\it lp}V \cdot F$ factorization is 
$\#{\mathcal P}$-hard, even when circuits for the rank functions of the 
domain code and image code are given.  

\smallskip

To define the class $\#{\mathcal P}$ we consider functions of the form 
$f: A^* \to \{0,1\}^*$, where $A$ is a finite alphabet, and elements of 
$\{0,1\}^*$ are interpreted as non-negative integers in binary representation.
Intuitively, for a function $f$ in $\#{\mathcal P}$ and for $x \in A^*$,  
$f(x)$ is the number of ways a relation that is parameterized by $x$ can
be satisfied. More precisely we will use the following definition of the 
$\#{\mathcal P}$; see e.g.\ \cite{Handb}.

\begin{defn} \label{numberP} \ 
A function $f: A^* \to \{0,1\}^*$ is in $\#{\mathcal P}$ \ iff \ there is a 
relation $R \subseteq A^* \times B^*$ (where $B$ is a finite alphabet) such 
that \\  
(1) for all $x \in A^*$: \ $f(x) \ = \ |\{ w \in B^* : (x,w) \in R\}|$, \ 
 with $f(x) \in \{0,1\}^*$ interpreted as an integer; \\ 
(2) $R$ is in ${\mathcal P}$ (deterministic polynomial time);  \\  
(3) $R$ is polynomially balanced (also called ``polynomially honest''); 
 i.e., there is a polynomial $p(.)$ such that for all $(x,w) \in R$, \    
$|x| \leq p(|w|)$ and $|w| \leq p(|x|)$.
\end{defn}

\begin{thm} \label{numberP_rank} {\bf (Ranking problem for finite maximal
prefix codes)} \\      
The following problem is $\#{\mathcal P}$-complete. \\  
{\sf Input:} A $\{0,1,\bot\}$-valued circuit that accepts a finite maximal 
prefix code $P \subset \{0,1\}^*$, and $x \in P$. \\  
{\sf Output:} The rank of $x$ in $P$ according to dictionary order.
\end{thm}
{\bf Proof.} The problem is clearly in $\#{\mathcal P}$ since 
${\rm rank}_P(x)$ is the number of words $w \in B^*$ (here $B = \{0,1\}$)
satisfying the relation ``$w \in P$ and $w <_d x$''. Moreover, the prefix 
code $P$ is given by a circuit, whose size is counted as part of the input 
size of the problem, so the relation ``$w \in P$ and $w <_d x$'' can be 
verified in deterministic polynomial time.

Next, we will reduce the $\#{\mathcal P}$-complete problem \#SAT to our 
problem.  For a boolean formula $\beta(x_1, \ldots, x_n)$ with $n$ boolean 
variables, let $T \subseteq \{0,1\}^n$ be the set of truth-value assignments 
that make $\beta$ true. Although $T$ is a finite prefix code, $T$ is not 
maximal, and the cardinality $|T|$ is not necessarily a power of 2; however, 
finding $|T|$, given $\beta$, is precisely the $\#{\mathcal P}$-complete 
problem \#SAT. We will use $T$ to construct a finite maximal prefix code $P$ 
(with $|P|$ a power of 2), whose ranking function determines $|T|$.
We use the notation ${\overline T} = \{0,1\}^n - T$.  Let \ 

\smallskip

$P_T \ = \ 00{\overline T} \ \cup \ 00T0 \ \cup \ 00T1 \ \cup \ $
           $01\, \{0,1\}^n \ \cup \ 1T \ \cup \ $
           $ 1 {\overline T}0 \ \cup \ 1 {\overline T}1$.          

\smallskip

\noindent Then $P_T$ is a finite maximal prefix code of cardinality
$|P_T| =  2^{n+2}$.  Membership in $P_T$ is easily decided by the formula
$\beta$.

Finally, $|T|$ is easily derived from the rank of $00 1^n$ or of
$001^n1$ in $P_T$. Indeed, if $1^n \in {\overline T}$ then \   
${\rm rank}_{P_T}(001^n) +1 = |T| \cdot 2 + |{\overline T}| = |T| + 2^n$; 
if $1^n \in T$ then \  ${\rm rank}_{P_T}(001^n1) +1 = |T| + 2^n$. 
Hence, $|T|$ can easily be obtained from ${\rm rank}_{P_T}(01^n)$ or 
${\rm rank}_{P_T}(01^n1)$; the numbers are written in binary, so the
representation of $2^n$ is not large.
 \ \ \ $\Box$

\begin{thm} \label{numberP_fact1} {\bf (${\it lp}V \cdot F$ factorization
problem, given $\Phi: P \to \{0,1\}^n$)} \\    
The following problem is $\#{\mathcal P}$-complete. \\
{\sf Input:} A $\{0,1,\bot\}$-valued circuit that computes a bijection 
$\Phi: P \to \{0,1\}^n$ (where $P$ is a finite maximal prefix code over 
$\{0,1\}$), and $x \in P$. \\  
{\sf Output:} The rank of $x$ in $P$ according to dictionary order. (Recall
that ${\rm rank}_P(.)$ is the $F$-part in the ${\it lp}V \cdot F$ 
factorization of $\Phi$.) 

\smallskip

The problem remains $\#{\mathcal P}$-complete if we assume that circuits
for both $\Phi$ and $\Phi^{-1}$ are given.
Also, evaluating $\pi$ or $\pi^{-1}$ is $\#{\mathcal P}$-complete
(where $\Phi(.) = \pi \, f(.)$ is the ${\it lp}V \cdot F$ factorization).
\end{thm}
{\bf Proof:} The problem is in $\#{\mathcal P}$ because the circuit for 
$\Phi$ can also be used to test membership in $P$.    
To show $\#{\mathcal P}$-hardness, let $P_T$ be as in Theorem 
\ref{numberP_rank} above, where $T \subseteq \{0,1\}^n$ is the set of truth
value assignments that make a given boolean formula $\beta$ true; again,
${\overline T}$ denotes $\{0,1\}^n - T$. 

\medskip

$P_T \ = \ 00{\overline T} \ \cup \ 00T0 \ \cup \ 00T1 \ \cup \ $
           $01\, \{0,1\}^n \ \cup \ 1T \ \cup \ $
           $ 1 {\overline T}0 \ \cup \ 1 {\overline T}1$.

\medskip

\noindent Let $\Phi: P_T \to \{0,1\}^{n+2}$ be the bijection defined as 
follows for all $x \in \{0,1\}^n$:

\smallskip

$00x \in 00{\overline T} \ \longmapsto \ 11x \in 11{\overline T}$

\smallskip

$00x0 \in 00T0 \ \longmapsto \ 0x0 \in 0T0$

\smallskip

$00x1 \in 00T1 \ \longmapsto \ 0x1 \in 0T1$

\smallskip

$01x \in 01 \, \{0,1\}^n \ \longmapsto \ 10x \in 10 \, \{0,1\}^n$

\smallskip

$1x0 \in 1 {\overline T}0 \ \longmapsto \ 0x0 \in 0 {\overline T}0$

\smallskip

$1x1 \in 1 {\overline T}1 \ \longmapsto \ 0x1 \in 0{\overline T}1$

\smallskip

$1x \in 1 T \ \longmapsto \ 11x \in 11 T$.

\smallskip

\noindent Clearly $\Phi$ and $\Phi^{-1}$ can easily be computed from the 
boolean formula $\beta$, and they have small circuits that can be derived
from the boolean formula $\beta$. 

Let $\Phi = \pi \, f$ be the ${\it lp}V \cdot F$ factorization of 
$\Phi$; then $f = {\rm rank}_{P_T}(.)$. We saw in Theorem
\ref{numberP_rank} above that evaluating ${\rm rank}_{P_T}(.)$ is a 
$\#{\mathcal P}$-complete problem.  Thus by the reduction of $f$ to $f^{-1}$
in Theorem \ref{inverting_f}, the problem of computing $f^{-1}$ is 
also $\#{\mathcal P}$-complete. 

To show that the evaluations of $\pi$ and $\pi^{-1}$ are 
$\#{\mathcal P}$-hard, note that $f = \pi^{-1} \Phi$ and 
$f^{-1} = \Phi^{-1} \pi$; since $\Phi$ and $\Phi^{-1}$ are easy to
evaluate, this reduces the $\#{\mathcal P}$-complete evaluation problems for 
$f$ and $f^{-1}$ to the evaluation of $\pi^{-1}$, respectively $\pi$.  
 \ \ \ $\Box$

\medskip

\noindent
The above Theorem means that ranking in $P_T$ according to the dictionary 
order is hard, but there may exist another bijection $P_T \to \{0,1\}^n$, 
namely $\Phi$, which provides an easy ranking in $P_T$.

\begin{thm} \label{numberP_fact2} {\bf (${\it lp}V \cdot F$ factorization,
given $\varphi: P_0 \to Q_0$, rank$_{P_0}$(.) and rank$_{Q_0}$(.))} \\
The following problem is $\#{\mathcal P}$-complete. \\
{\sf Input},  consisting of three parts: \\ 
$\bullet$ A $\{0,1,\bot\}$-valued circuit that computes a bijection
$\varphi: P_0 \to Q_0$ (where $P_0$ and $Q_0$ are finite maximal prefix codes
over $\{0,1\}$), \\  
$\bullet$ two $\{0,1,\bot\}$-valued circuits that compute the rank functions
of $P_0$, respectively  $Q_0$, \\ 
$\bullet$ and $x \in P_1$ (where $P_1$ is the domain code that $\varphi$ 
receives when it is restricted so as to have {\it imC}$(\varphi)$
$= \{0,1\}^n$, where $n = {\rm max}\{ |q| : q \in Q_0\}$). \\  
{\sf Output:} The rank of $x$ in $P_1$ according to dictionary order.

\smallskip

\noindent
The problems remains $\#{\mathcal P}$-complete if we assume that circuits
for both $\varphi$ and $\varphi^{-1}$ are given.

Also, evaluating $\pi$ or $\pi^{-1}$ is $\#{\mathcal P}$-complete
(where $\varphi = \pi \, f$ be the ${\it lp}V \cdot F$ factorization).
\end{thm}
{\bf Proof:} The problem is in $\#{\mathcal P}$ because the circuit for
$\varphi$ can also be used to obtain a circuit for the restriction
$P_1 \to \{0,1\}^n$ of $\varphi$, by Lemma \ref{makingImC_n}; this circuit can
then be used to test membership in $P_1$.

To show $\#{\mathcal P}$-hardness, let $T$ and $P_T$ be as in the proofs 
of Theorems \ref{numberP_rank} and \ref{numberP_fact2}. Let
 
\smallskip

$P_0 \ = \ \{00,01,1\} \cdot \{0,1\}^n$
    $ \ = \ 00T \ \cup \ 00{\overline T} \ \cup \ 01 \, \{0,1\}^n $
       $ \ \cup \ 1T \ \cup \ 1{\overline T}$,  

\smallskip

$Q_0 \ = \ \{0,10,11\} \cdot \{0,1\}^n$
    $ \ = \ 0T \ \cup \ 0{\overline T} \ \cup \ 10 \, \{0,1\}^n $
       $ \ \cup \ 11T \ \cup \ 11{\overline T}$,

\smallskip

\noindent and define $\varphi: P_0 \to Q_0$ by 

\smallskip

$00x \in 00{\overline T} \ \longmapsto \ 11x \in 11{\overline T}$

\smallskip

$00x \in 00T \ \longmapsto \ 0x \in 0T$ 

\smallskip

$01x \in 01\, \{0,1\}^n \ \longmapsto \ 10x \in 10 \, \{0,1\}^n$

\smallskip

$1x \in 1{\overline T} \ \longmapsto \ 0x \in 0{\overline T}$.

\smallskip

$1x \in 1T \ \longmapsto \ 11x \in 11T$

\smallskip

\noindent Then $\varphi$, ${\rm rank}_{P_0}(.)$, and ${\rm rank}_{Q_0}(.)$,
and their inverses have small circuits, that are easily derived from the
boolean formula $\beta$. 

Next, we restrict $\varphi$ in such a way that its image code becomes
$\{0,1\}^{n+2}$. The resulting bijection is exactly the bijection
$\Phi: P_T \to \{0,1\}^{n+2}$ of the proof of Theorem \ref{numberP_fact1}.
All the claimed conclusions of Theorem \ref{numberP_fact2} now follow
from Theorem \ref{numberP_fact1}. 
 \ \ \ $\Box$

\bigskip

\noindent The above $\#{\mathcal P}$-completeness results imply that finding
circuits for the ${\it lp}V \cdot F$ factors $\pi, f$ of $\varphi \in V$ is 
difficult (if ${\mathcal P} \neq NP$, etc.). However, whether this implies 
that the factors require large circuits remains a very difficult open problem.


\bigskip

\bigskip



\bigskip

\bigskip

\noindent {\bf Jean-Camille Birget} \\
Dept.\ of Computer Science \\
Rutgers University at Camden \\
Camden, NJ 08102, USA \\
{\tt birget@camden.rutgers.edu}

\end{document}